\documentclass[12pt]{article}
\usepackage{epsfig,amssymb,epic,eepic,color, graphics,amsmath}
\usepackage[english]{babel}
\input{amssym.def}

\def\nd{\noindent}

\newenvironment{demo}{\nd {\bf Proof: }}{${}$\hfill $\diamond$ \medskip}

\newtheorem{theo}{Theorem}

\newtheorem{defi}{Definition}

\newtheorem{stat}{Statement}

\begin{document}
\sloppy
\date{\today}
\title{On 2-diffeomorphisms with one-dimensional  basic sets and a finite number of moduli}
\author{V.\,Z.~Grines\footnote{National Research University Higher School of Economics,  vgrines@yandex.ru}, O.\,V.~Pochinka\footnote{National Research University Higher School of Economics, olga-pochinka@yandex.ru}, S.\,Van~Strien\footnote{Imperial College, 	
s.van-strien@imperial.ac.uk}}

\maketitle

\begin{abstract} This paper is a step towards the complete topological classification 
of $\Omega$-stable diffeomorphisms on an orientable closed surface, aiming to 
give necessary and sufficient conditions for two such diffeomorphisms to be topologically
conjugate without assuming that the diffeomorphisms are necessarily close to each other.
In this paper we will establish such a classification within a certain class $\Psi$ of $\Omega$-stable
diffeomorphisms defined below. To determine whether two diffeomorphisms from this 
class $\Psi$ are topologically conjugate, we give 
(i) an algebraic description of the dynamics on their non-trivial basic sets,  (ii)
a geometric description of how invariant manifolds intersect, and (iii) 
define numerical invariants, called moduli, associated  to 
orbits of tangency of stable and unstable manifolds of saddle periodic orbits. 
This description determines the {\em scheme} of a diffeomorphism,
and we will show that two diffeomorphisms from $\Psi$ are topologically conjugate
if and only if their schemes agree. 
\end{abstract}



Bibliography: \ref{last} names.

\section{Introduction and formulation of the results}
The topological classification of structurally stable diffeomorphisms on closed orientable surfaces has made tremendous progress in the last 25 years, due to the work of C. Bonatti,  V. Grines, R. Langevin, A. Zhirov, R. Plykin and etc. (see for example \cite{ArGr90}, \cite{BoLa}, \cite{BoGrLa01}, \cite{GrZhu2006} for the history of the subject and more information). Any such classification naturally includes a description of its basic sets and a non-trivial description of how invariant manifolds of periodic points intersect. If invariant manifolds of saddles have tangencies, then the topological classification also involves expressions, called {\em moduli}, related to eigenvalues at saddle points, as was discovered by J. Palis \cite{Pa78}.

A first important step in the direction of a topological classification of $\Omega$-stable diffeomorphisms on orientable closed surfaces 
was made  by W. de Melo and S.J. van Strien in \cite{MeSt1987}, where they found necessary and sufficient conditions for $\Omega$-stable diffeomorphisms 
to have a  {\em finite number of moduli}. (A diffeomorphism $f$  is said to have a finite number of moduli if one can parametrise topological conjugacy classes of a neighbourhood of  $f$ by a finite number of parameters). Their result is local in the sense that it only considers 
the topological conjugacy of  two diffeomorphisms which are sufficiently close to each other. To deal with the global situation, T. Mitryakova and O. Pochinka \cite{MiPo2010} partly generalised the previous result by construction a complete invariant for $\Omega$-stable diffeomorphisms for a certain class of $\Omega$-stable diffeomorphisms (with at most a finite number of periodic points)  which can in general be ``far'' from each other.

Here we present the topological classification considering a wider class than 
in \cite{MiPo2010}, within this class  the existence of one-dimensional attractors and repellers is allowed. This class $\Psi$ will be defined formally below. 

Let $M^2$ be an orientable closed surface and $f\colon M^2\to M^2$ be an A-diffeomorphism, i.e. an Axiom A diffeomorphism. By S. Smale \cite{Smale67}, the non-wandering set $ NW(f)$ of $f$ is represented as a finite union of disjoint closed invariant sets $\Lambda_1, \ldots, \Lambda_k $, called {\it basic sets}, each of which contains a dense orbit. A basic set which consists of a periodic orbit will be called {\it trivial} and otherwise it is called {\it non-trivial}.  

Let $\Lambda$ be a one-dimensional basic set of $f$. By R. Plykin \cite{Plyk74}, $\Lambda$ is either an attractor or a repeller.  

\begin{figure}[ht] \centerline
{\epsfig{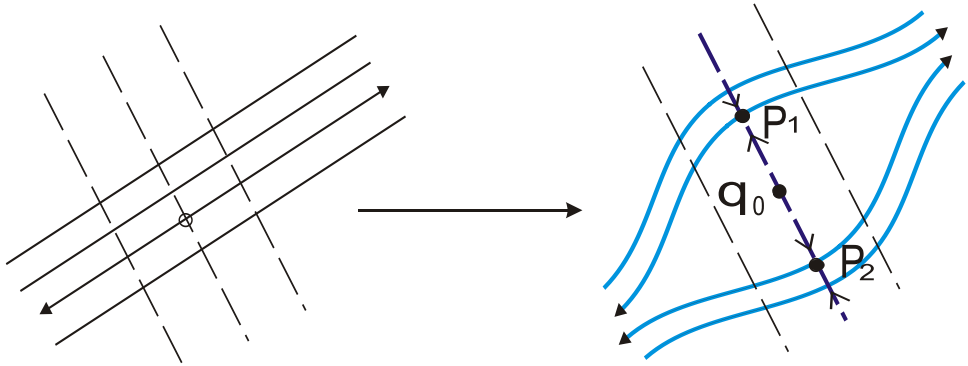}}\caption{DA-diffeomorphism} \label{DA}
\end{figure}

According to \cite{Gr74} (Definition 3) a point $p$ is called an {\it $s$-boundary ($u$-boundary) point} 
 of attractor (repeller) $\Lambda$, if one of the connected components of the set $W^s_p\setminus p\,(W^u_p\setminus p)$ is disjoint from $\Lambda$; denote by $\ell_p$ such a component (see Figure \ref{DA}, where the construction of a DA-diffeomorphism is represented and where $p_1,p_2$ are the $s$-boundary points).  

\begin{figure}[ht] \centerline
{\epsfig{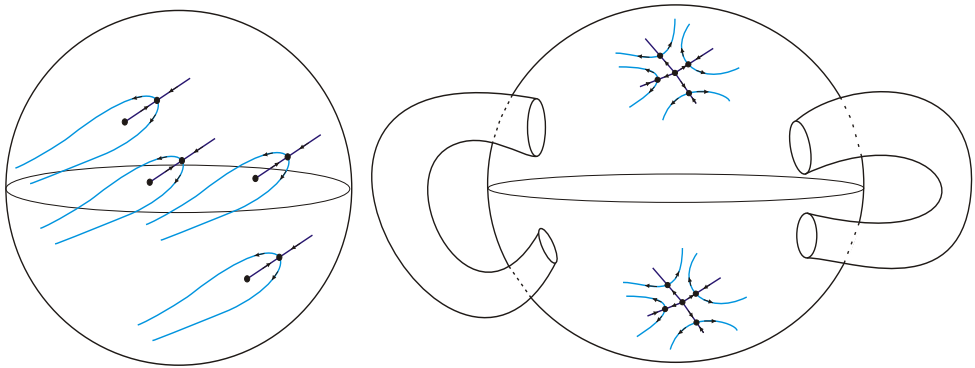}}
\caption{Other examples of 1-dimensional attractors} \label{Pl}
\end{figure}

For definiteness we suppose that $\Lambda$ is an attractor (all notions for repeller can be obtained by considering $f^{-1}$). By \cite{Gr75} (Lemmas 2.4, 2.5), each $s$-boundary point is necessarily periodic and the set $\Lambda$ has a non-empty and finite set of $s$-boundary points\footnote{In fact the existence and finiteness of the set of boundary points without the term ``boundary point'' was proved by S. Newhouse and J. Palis in \cite{NP73} (Proposition 1).}. We denote this set by $P_\Lambda$. 

\begin{defi}[Separable one-dimensional attractors]  \label{ott} We say that a 1-dimensional attractor $\Lambda$ of an $A$-diffeomorphism $f$ is {\em separable} if a union $Y_\Lambda$ of saddle and source trivial basic sets of the diffeomorphism $f$ exists with the following properties:  

1) $cl(W^s_\Lambda)\setminus W^s_\Lambda=W^s_{Y_\Lambda}$;

2)  $cl(\ell_p)\setminus \ell_p=p\cup\alpha$ for every $s$-boundary point $p\in P_\Lambda$,  where $\alpha\in Y_\Lambda$ is a source point;

3) for every saddle point $\sigma\in Y_\Lambda$ the manifold $W^s_\sigma$ does not contain  heteroclinic points. 
\end{defi}

\begin{figure}[ht] \centerline
{\epsfig{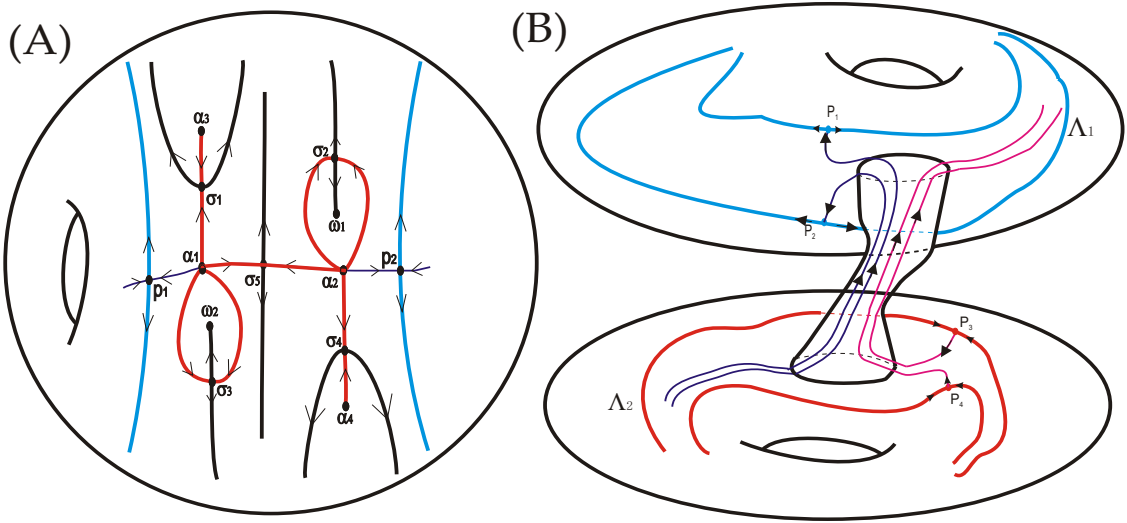}}\caption{\small In these figures $\alpha,\omega, \sigma$ denote periodic points of source, sink and saddle type. In Figure (A) a separable one-dimensional attractor on the torus is shown, where $cl~(W^s_\Lambda)\setminus W^s_\Lambda$
is the region bounded by the dark (red) curves and $\ell_{p_i}$ are the curves connecting
$p_i$ to $\alpha_i$; In Figure (B) a situation is shown with 
a non-seperable one-dimensional attractor and repeller on the pretzel; here the curves $\ell_{p_i}$  do not land on repelling fixed points, and so condition 2) in the definition of separable one-dimensional attractor is violated.} 
\label{kart}
\end{figure}

This definition is illustrated in Figure \ref{kart}.

It follows from \cite{Gr97} (Lemma 1, Lemma 2) that any one-dimensional basic set of a structural stable diffeomorphism $f:M^2\to M^2$ is separable. We prove the following stronger result.  

\begin{theo}\label{1-sep} If an $\Omega$-stable diffeomorphism $f:M^2\to M^2$ has a finite number of moduli then any of its one-dimensional basic set is separable.  
\end{theo}

The proof of Theorem \ref{1-sep} is based on  necessary and sufficient conditions, found in \cite{MeSt1987},  under which a diffeomorphism of an orientable surface has a finite number of moduli of topological conjugacy, and described the structure of the neighborhood of such a diffeomorphism. 

\begin{stat} [Criteria of a finite number moduli, \cite{MeSt1987}] \label{fimo} Let $f:M^2\to M^2$ be  an  $\Omega$-stable $C^2$-diffeomorphism. Then $f$ has a finite number moduli if and only if it satisfies the conditions below:

(1) if $x,y\in NW(f)$ are such that $W^u_x$ is not transverse to $W^s_y$ then the basic 
sets containing $x$ and $y$ are trivial;

(2) there is only a finite number of orbits of non-transverse intersections between 
stable and unstable manifolds and the contact between these manifolds along each 
of these orbits is of finite order;

(3) if $p,q$ are periodic points from trivial basic sets  such that $W^u_p$ has an orbit of non-transverse intersection with $W^s_q$ then the number of orbits in $W^s_p$ (resp. in $W^u_q$) belonging to some unstable (resp. stable) manifolds of periodic saddle points of $f$ is finite;

(4) if $x$ is a point of non-transverse intersection of $W^u_p$ and $W^s_q$ then there 
exists an arc $\Sigma$ transverse to $W^u_p$ at $x$ such that no connected component of
$\Sigma\setminus \{x\}$ contains points of both stable and unstable manifolds of saddles;

(5) if $W^u_p$ has a point of non-transverse intersection with $W^s_q$,and $W^u_q$ 
has a point of non-transverse intersection with $W^s_r$, then there is no saddle point
of $f$ whose unstable manifold (resp. stable manifold) intersects $W^s_p$ (resp. $W^u_r$).
\end{stat}

\begin{defi}[The class $\Psi$]
An  orientation preserving $\Omega$-stable $C^2$-diffeomorphism $f:M^2\to M^2$
is called a diffeomorphism of class $\Psi$ if it has a finite number of moduli and the following properties are satisfied:

1) each non-trivial basic set $\Lambda$ of $f$ is one-dimensional;

2) heteroclinic orbits can be contained in the stable or the unstable manifold of a periodic point of the trivial basic set of $f$, but not in both. 
\end{defi}

\begin{figure}[ht] \centerline
{\epsfig{file=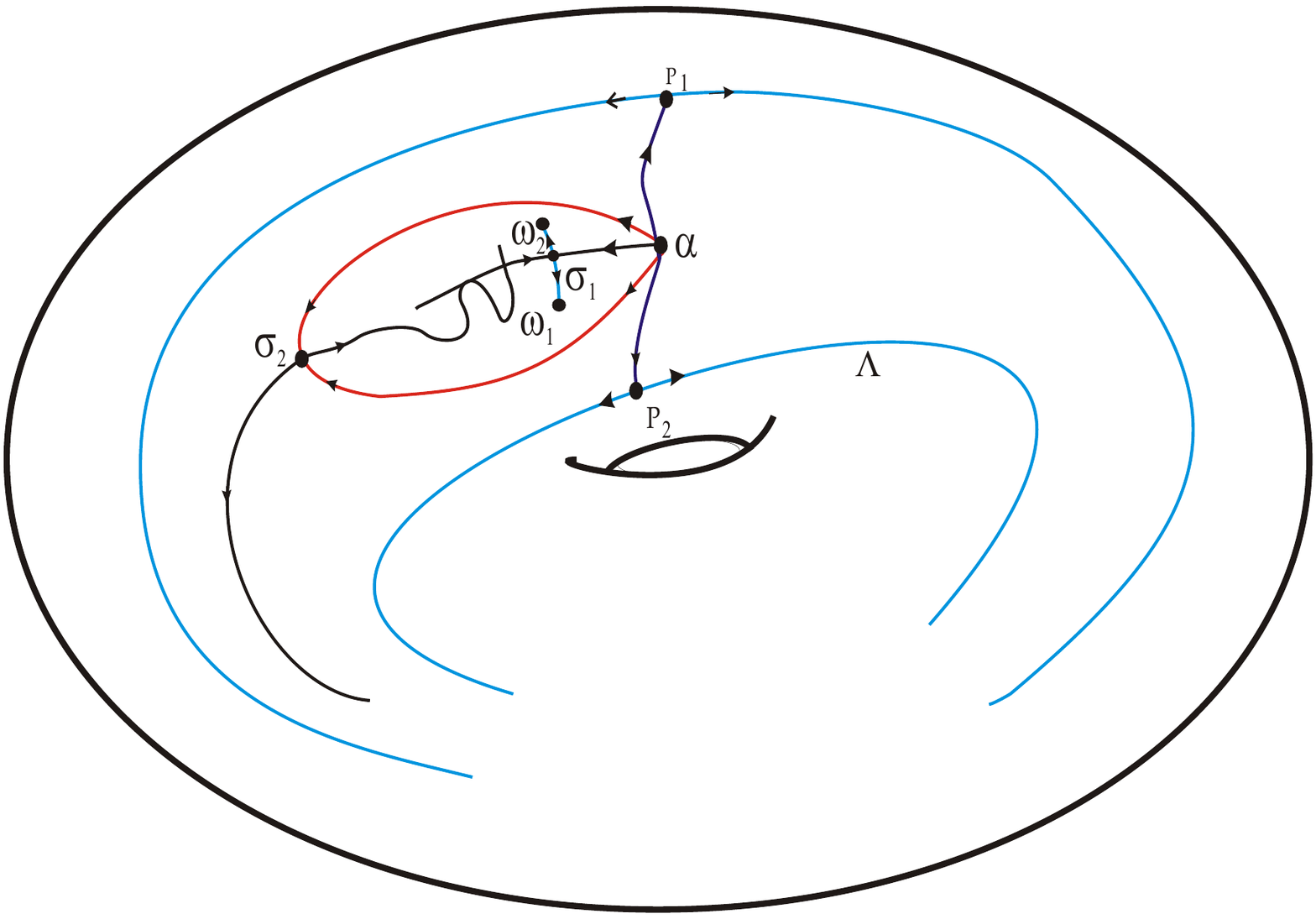,width=10cm,height=5.5cm}}
\caption{An example of a diffeomorphism $f$ from $\Psi$} \label{UL}
\end{figure}

Let $f\in\Psi$ and $W^u_x$ be not transverse to $W^s_y$ for some saddle periodic points $x,y\in NW(f)$. Set $\Theta_{xy}=\frac{ln|\lambda_x|}{ln|\mu_y|}$, where $\lambda_x$ is the eigenvalue of $D_f$ at $x$ which is less than one by absolute value and $\mu_y$ is the eigenvalue of $Df$ at $y$ which is  greater than one  by absolute value. Denote by $\Psi^*$ the set of diffeomorphisms $f\in\Psi$ such that $\Theta_{xy}$ is an irrational for any such pair $x,y$.  

In section \ref{BP} we introduce the notion of a scheme of diffeomorphism $f$ containing 

(i) an algebraic description of the dynamics on its non-trivial basic sets,  

(ii) a geometric description of how invariant manifolds intersect, 

(iii) numerical invariants, called moduli, associated  to 
orbits of tangency of stable and unstable manifolds of saddle periodic orbits

and define an equivalence of two schemes.

The main result of this paper is the following theorem.

\begin{theo}[Classification with class $\Psi$] $ $

1. If the schemes of diffeomorphisms $f,f'\in\Psi$ are equivalent, then the diffeomorphisms are topologically conjugate.

2. Diffeomorphisms $f,f'\in\Psi^*$ are topologically conjugate if and only if their schemes are equivalent.\label{th1}
\end{theo}

{\it Acknowledgements.} This work was supported by the Russian Foundation for Basic Research (project nos. 15-01-03687-a, 16-51-10005-Ko\_a), Russian Science Foundation  (project no 14-41-00044), the  Basic Research Program at the HSE (project 98) in 2016 and the European Union ERC AdG grant no 339523 RGDD.

\section{Descriptions of diffeomorphisms from $\Psi$}  
\label{BP}

\subsection{An algebraic description of the dynamics on one-dimensional basic set}
\label{Lam}

Now let $\Lambda$ be a 1-dimensional attractor of an A-diffeomorphism $f:M^2\to M^2$. From \cite{Bow71} $ \Lambda$ is represented as a finite union of disjoint compact sets $ \Lambda _ {1},\dots, \Lambda _ {k} $, which are cyclically transformed into each other under the action of f.
 Moreover, $cl~(W^s_x\cap\Lambda_i)=\Lambda_i$ and $cl~(W^u_x\cap\Lambda_i)=\Lambda_i$ for  any point $ x \in \Lambda _ {i} $. Every $ \Lambda _ {i} $ is called {\it periodic (or $ C $-dense) component} of the basic set $\Lambda$. In this section we suppose that the attractor $\Lambda$ consists of only one periodic component. We will now associate to $\Lambda$ a closed neighbourhood $N_\Lambda$. 

\begin{defi}[The bunch of an attractor] A {\em bunch $b$} of an attractor $\Lambda$ is the union of the maximal number $r_b$ of the unstable manifolds $W^{u}_{p_{1}},\dots,W^{u}_{p_{r_b}}$ of the $s$-boundary points $p_{1}, \dots, p_{r_b}$ of the set $\Lambda$ whose separatrices\footnote{{\it Stable (unstable) separatrix} of a hyperbolic periodic point $p$ is a connected component of the set $W^{s}_{p}\setminus p\,(W^{u}_{p}\setminus p)$.} $\ell_{p_1},\dots,\ell_{p_{r_b}}$ belong to the same connected component of the set $W^s_\Lambda\setminus\Lambda$. The number $r_b$ is said to be {\it the degree of the bunch}.
\end{defi}

Let $\delta\in\{u,s\}$ and $x\in\Lambda$. For points $y,z \in W^{{\delta}}_{x},~(y\neq z)$ let  $$[y,z]^{{\delta}},~~[y, z)^{{\delta}},~~( y, z]^{{\delta}},~~( y, z)^{{\delta}}$$ denote the connected arcs on the manifold $W^{{\delta}}_{ x}$ with the boundary points
$ y,z$.

Denote by $B_\Lambda$ the set of all bunches of $\Lambda$. From the definition of a bunch  $b\in B_\Lambda$ of degree $r_b$ it follows that there is a sequence of points $x_{1}, \dots,x_{2r_b}$ such that:

1) $x_{2j-1}, x_{2j}$ belong to the different connected components of the set
$W^{u}_{p_{j}} \setminus p_{j}$;

2) $x_{2j+1} \in W^{s}_{x_{2j}}$ (we set
$x_{2r_b+1} = x_{1}$)

3) $(x_{2j}, x_{2j+1})^{s} \cap \Lambda =
\emptyset$, $j=1,\dots,r_b$.

For each $j\in\{1,\dots,r_b\}$ we pick a pair of points $\tilde x_{2j-1},  \tilde  x_{2j}$ and a simple curve $\ell_{j}$  with  boundary points $\tilde x_{2j-1}, \tilde x_{2j}$ such that:

1) $(\tilde x_{2j},\tilde x_{2j+1})^{s}
\subset(x_{2j}, x_{2j+1})^{s}~$
$(x_{2r+1} = x_{1})$;

2) the curve   $\ell_{j}$  transversally intersects at a unique point the stable manifold of any point on the arc $(x_{2j-1}, x_{2j})^{u}$;

\begin{figure} \centerline
{\epsfig{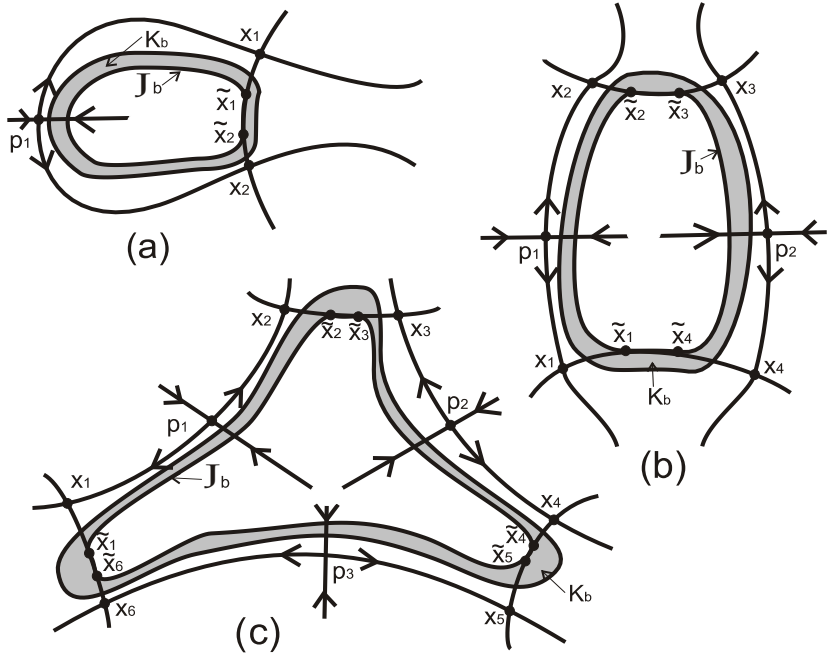}}
\caption{Construction of the surface $N_\Lambda$.}
\label{M-Lambda}
\end{figure}

3) ${L}_b=\bigcup\limits_{j=1}^{r_b}\left[\ell_{j} \cup (\tilde x_{2j},\tilde x_{2j+1})^{s}\right]$ is a simple closed smooth curve and the set ${L}_{\Lambda}=\bigcup\limits_{b\in B_\Lambda}{L}_b$ is such that: 

a) $f({L}_{\Lambda})
\cap {L}_{\Lambda}=\emptyset$;

b) for every curve $L_b,~b\in B_\Lambda$ there is a curve in the set $f({L}_{\Lambda})$ such that these two curves are the boundaries of an annulus $K_b$; 

c) the annuli  $\{K_b,~b\in B_\Lambda\}$ are pairwise disjoint (see Figure \ref{M-Lambda}).

Let  $N_\Lambda=\Lambda\cup\bigcup\limits_{n\geq 1}f^n(\bigcup\limits_{b\in B_\Lambda}K_b)$. By construction the annuli $\{K_b,~b\in B_\Lambda\}$ consist of the wandering points of the diffeomorphism $f$,  $N_\Lambda$ is a surface with non-empty boundary and $N_\Lambda$ is a neighbourhood of the attractor $\Lambda$, which we call the support of $N_\Lambda$.

Let $p_{_{\Lambda}}:\mathbb U_{N_{\Lambda}}\to N_\Lambda$ be the universal covering where $\mathbb U_{N_{\Lambda}}$ is a subset of Lobachevsky plane and let $\mathbb G_{N_{\Lambda}}$ be the group of its covering transformations. Set $\mathbb E_{N_\Lambda}=\partial \mathbb U_{N_{\Lambda}}$. A lift $\bar f_{\Lambda}:\mathbb U_{N_{\Lambda}}\to\mathbb U_{N_{\Lambda}}$ of $f_\Lambda=f|_{N_\Lambda}$ with respect to  $p_{_{\Lambda}}$ induces an automorphism $T_{\bar f_{\Lambda}}$ of the group $\mathbb G_{N_{\Lambda}}$ which acts by the formula $T_{\bar f_{\Lambda}}(g)=\bar f_{\Lambda}g\bar f_{\Lambda}^{-1}$. Set $\bar\Lambda=p^{-1}_{\Lambda}(\Lambda)$. If $x\in\Lambda $ then let   $\bar x\in\bar\Lambda $ denote a point in the preimage $p_{_{\Lambda}}^{-1}(x)$ and let $w_{\bar x}^\delta$ be the connected component of $p^{-1}(W_x^{\delta})$ containing $\bar x$.  Let us choose a parametrisation  $\mathbb R\ni t\to W_x^{\delta}(t)$ of $W_x^{\delta}$   such that $W_{x}^{\delta}(0)=x$. Then
$w_{\bar x}^{\delta}(t)$ is a point on $w_{\bar x}^{\delta}$ 
such that $p_{_{\Lambda}}(w_{\bar x}^{\delta}(t))=W_{x}^{\delta}(t)$ and $W_{x}^{{\delta}+},~W_{x}^{{\delta}-}~(w_{\bar x}^{{\delta}+},~w_{\bar x}^{{\delta}-})$ are the connected components of the curve $W_{x}^{{\delta}}\setminus x~(w_{\bar x}^{{\delta}}\setminus \bar x)$ for $t>0$, $t<0$ respectively. 
For points $\bar y,\bar z \in w^{{\delta}}_{\bar x},~(\bar y\neq\bar z)$ let $[\bar y,\bar z]^{{\delta}}$, $[\bar y,\bar z)^{{\delta}}$, $(\bar y,\bar z]^{{\delta}}$,  $(\bar y,\bar z)^{{\delta}}$ denote the connected arcs on the manifold
$w^{{\delta}}_{\bar x}$ with boundary points
$\bar y,\bar z$.

\begin{defi} [Asymptotic direction] We say that a curve $w_{\bar x}^{{\delta}{\nu}}$ has the asymptotic direction ${\delta}_{\bar x}^{{\nu}}$ for $t\to {\nu}\infty,~\nu\in\{-,+\}$ if  $cl(w_{\bar x}^{{\delta}{\nu}})\setminus w_{\bar x}^{{\delta}{\nu}}$ is equal to 
$\bar x$ and ${\delta}_{\bar x}^{{\nu}}\in \mathbb E_{N_\Lambda}=\partial \mathbb U_{N_\Lambda}$.
\end{defi}

As before, let $P_\Lambda$ be the set of $s$-boundary points of $\Lambda$.
For a boundary point $p\in P_\Lambda$ denote by $\ell_{p}^{\infty}$ a connected component of $W^s_p\setminus p$ different from $\ell_p$ and by $\ell_{\bar p}^{\infty}$ the connected component of $w^s_{\bar p}\setminus \bar p$ for which $p_{_\Lambda}(\ell_{\bar p}^{\infty})=\ell_{p}^{\infty}$. 

\medskip 

It was proved in \cite{Gr1999} that for each point $x\in(\Lambda\setminus W^s_{P_\Lambda})$ the curve $w_{\bar x}^{s}$ has two distinct boundary points (asymptotic directions) $s_{\bar x}^{+},~s_{\bar x}^{-}$. Finally, for every point  $p\in {P_\Lambda}$ the curve $w_{\bar p}^{s\infty}$ has an asymptotic direction $s_{\bar p}^{\infty}$.

\begin{stat} [Conjugacy on one-dimensional attractors, \cite{Gr1999}] Let ${\Lambda}$, ${\Lambda'}$ be attractors such that there is an automorphism $\psi_{_\Lambda}:\mathbb G_{{N_\Lambda}}\to\mathbb G_{N_{\Lambda'}}$ with property  $T_{{\bar f'_{\Lambda'}}} =
\psi_{_\Lambda} T_{{\bar f_\Lambda}}\psi^{-1}_\Lambda$. Then 
 
1) $\psi_{_\Lambda}$ is uniquely induces a homeomorphism $\psi^*_\Lambda:\mathbb E_{N_\Lambda}\to\mathbb E_{N_{\Lambda'}}$;  

2) for every point $\bar x\in\bar\Lambda$ there is a unique point $\bar x'\in\bar\Lambda'$ such that $\psi^*_\Lambda(cl~w^\delta_{\bar x}\cap\mathbb E_{N_\Lambda})=cl~w^\delta_{\bar x'}\cap\mathbb E_{N_{\Lambda'}}$ for $\delta\in\{u,s\}$ and the map $\bar{\varphi}_{\Lambda}:\bar\Lambda\to\bar\Lambda'$, assigning $\bar x'$ to $\bar x$, is a homeomorphism;

3) $\bar{\varphi}_{\Lambda}$ induces the homeomorphism $${\varphi}_{\Lambda}=p_{_{\Lambda'}}\bar{\varphi}_{\Lambda} p^{-1}_{_{\Lambda}}:\Lambda\to\Lambda'$$ conjugating $f|_\Lambda$ with $f'|_{\Lambda'}$ and possesses the property: if $a,b\in W^s_x,~x\in\Lambda$ then ${\varphi}_{\Lambda}(a),{\varphi}_{\Lambda}(b)\in W^s_{{\varphi}_{\Lambda}(x)}$.

\label{atr.l8}
\end{stat}

It immediately follows from Statement \ref{atr.l8} that each isomorphism $\psi_{_\Lambda}$ with property  $T_{{\bar f'_{\Lambda'}}} =
\psi_{_\Lambda} T_{{\bar f_\Lambda}}\psi^{-1}_\Lambda$ uniquely induces a one-to-one map $$\hat\psi_{_{\Lambda}}:P_\Lambda\to P_{\Lambda'}.$$

\subsection{Moduli associated to diffeomorphisms from the class $\Psi$}

For two diffeomorphisms from $\Psi$ to be topologically conjugate certain moduli conditions have to be satisfied. Let us
define these conditions now. For $f\in \Psi$ denote by $\Omega_f$ the set of trivial basic sets of $f$ and by $\Omega^0,\Omega^1,\Omega^2$ its subsets consisting of the sinks, saddles and sources, accordingly. 
For a saddle point $\sigma\in\Omega^1$ of a diffeomorphism $f\in \Psi$ denote by $k_\sigma$ the period of $\sigma$ and  $\mu_{\sigma}$,  $\lambda_{\sigma}$ denote the eigenvalue of $Df^{k_\sigma}_{\sigma}$ which  are greater and less than one by absolute value, accordingly ($|\mu_{\sigma}|>1>|\lambda_{\sigma}|>0$). 

For $0<|\lambda|<1<|\mu|$ denote by $f_{\mu,\lambda}:\mathbb{R}^2\to\mathbb{R}^2$ linear diffeomorphism given by the formula $$f_{\mu,\lambda}(x,y)=(\mu x,\lambda y).$$ Set $$U_{\mu,\lambda}=\{(x,y)\in\mathbb{R}^2: \vert x\vert {\vert y\vert}^{-\log_{_{\lambda}}\mu} \leq 1\}.$$ Notice that the set $U_{\mu,\lambda}$ is $f_{\mu,\lambda}$-invariant and possesses two $f_{\mu,\lambda}$-invariant foliations $\mathcal F^s=\bigcup\limits_{c\in\mathbb R}\{(x,y)\in U_{\mu,\lambda}:x=c\}$ and $\mathcal F^u=\bigcup\limits_{c\in\mathbb R}\{(x,y)\in U_{\mu,\lambda}:y=c\}$. 

\begin{figure}[ht] \begin{center}\epsfig
{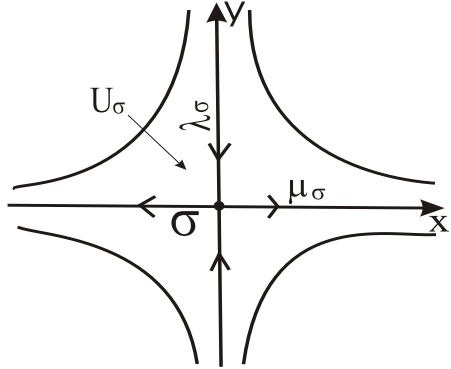} \caption{Linearizable neighborhood $U_{\sigma}$} \label{neib}
\end{center}
\end{figure}

\begin{defi}[$C^1$ linearization]  A saddle point $\sigma\in\Omega^1$  and an $f^{k_\sigma}$-invariant neighbourhood  $U_{\sigma}$ of $\sigma$ form a $C^1$-linearization (see Figure \ref{neib}) if:

1) there is a $C^1$-diffeomorphism $\psi_{\sigma}:U_{\sigma}\to U_{\mu_{\sigma}, \lambda_{\sigma}}$ conjugating $f^{k_\sigma}\vert_{U_{\sigma}}$ with $f_{\mu_{\sigma},\lambda_{\sigma}}|_{U_{\mu_{\sigma},\lambda_{\sigma}}}$;

2) each leaf of the foliations $\mathcal F^s_\sigma=\psi_\sigma^{-1}(\mathcal F^s),~\mathcal F^u_\sigma=\psi_\sigma^{-1}(\mathcal F^u)$ is $C^2$-smooth.
\label{lin}
\end{defi}

The existence of a linearizable neighborhood for any saddle point of a diffeomorphism $f$ from $\Psi$ (or indeed any $C^2$ diffeomorphism) is well-known, see for example \cite{ChShShTu} (Chapter 5).

Denote by $\mathcal{A}$ the set of points at which one-sided heteroclinic tangencies of invariant manifolds of saddle points of the diffeomorphism $f$ there is. For each $a\in \mathcal{A}$ denote by $\sigma^u_a\in\Omega^{1}$, $\sigma^s_a\in\Omega^{1}$ the saddle points such that $a\in W^s_{\sigma_a^s}\cap W^u_{\sigma_a^u}$. Set  $\mu_{a}=\mu_{\sigma_a^s}$ and $\lambda_{a}=\lambda_{\sigma_a^u}$. 

\begin{figure}[ht]
\begin{center}\epsfig
{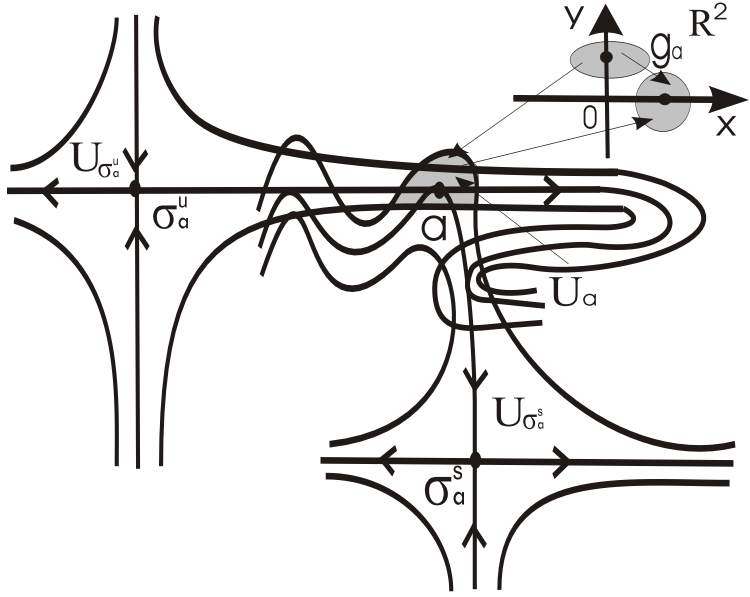} 
\caption{\small The transition map $g_a$}\label{map}
\end{center}
\end{figure}

Denote by $U_a$ the connected component of  $U_{\sigma_a^s}\cap
U_{\sigma_a^u}$ containing the point $a$. For any point $z \in U_{a}$ let us set  $z^s=\psi_{\sigma_{a}^s}(z)=(z^s_{x},z^s_{y})$ and $z^u=\psi_{\sigma_{a}^u}(z)=(z^u_{x},z^u_{y})$. Set $g_a=\psi_{\sigma_a^u}  ({\psi_{\sigma_a^s}\vert_{U_a}})^{-1}:\psi_{\sigma_a^s}(U_a)\to \psi_{\sigma_a^u}(U_a)$ (see Figure \ref{map}) and write the map $g_a$ in the coordinate form $$g_a(x,y)=(\xi_a(x,y), \eta_a(x,y)).$$ 

Set $$\tau_{a}=\frac{\partial{\eta}_{a}} {\partial{x}}(a^{s}).$$

Set $$H_a=\mathcal{A} \cap W^s_{\sigma_a^s}\cap W^u_{\sigma_a^u}.$$

\begin{stat} [Moduli for $f\in\Psi$, \cite{MeSt1987}, \cite{MiPo2010}] \label{mnogo} 
 $ $
 
\begin{enumerate}
\setlength\itemsep{0em}
\item Let $f\in \Psi$, $a \in \mathcal{A}$ and $k \in \mathbb{Z}$. Then $$\tau_{f^k(a)}=\left|\frac{\lambda_{a}}{\mu_{a}}\right|^k\cdot \tau_a.$$\label{S4}
\item If diffeomorphisms $f,f'\in \Psi$ are topologically conjugate by means of a homeomorphism $h$ such that $h(a)=a'$ for a point $a \in \mathcal{A}$ and $h(\sigma^s_a)=\sigma^s_{a'}$,
$h(\sigma^u_a)=\sigma^u_{a'}$ then $$\frac{\ln|\lambda_{a}|}{\ln|\mu_{a}|}=\frac{\ln|\lambda_{a'}|}{\ln|\mu_{a'}|}.$$ \label{S1}
\item  If diffeomorphisms $f,f'\in \Psi^*$ are topologically conjugate by means of a homeomorphism $h$ such that  $h(\sigma^s_a)=\sigma^s_{a'}$,  
$h(\sigma^u_a)=\sigma^u_{a'}$ for some points $a \in \mathcal{A},~a' \in \mathcal{A}'$ and $h(a_1)=a'_1,~h(a_2)=a'_2$ for some points $a_1,a_2\in H_a,~a'_1,a'_2\in H_{a'}$ then  $$\left|\frac{\tau_{a_2}}{\tau_{a_1}}\right|^{\frac{1}{\ln|\mu_{a}|}}=\left|\frac{\tau_{a'_2}}{\tau_{a'_1}}\right|^{\frac{1}{\ln|\mu_{a'}|}}.$$ 
 \label{S3}
\end{enumerate}
\end{stat}

\subsection{Geometric description of the intersection pattern of
invariant manifolds} 

Denote by $\mathcal L^s,\mathcal L^u$ the sets of non-trivial attractors respectively repellers. Set $\mathcal L=\mathcal L^s\cup\mathcal L^u$. As before let $\Omega^0,\Omega^1,\Omega^2$ be the sets of sinks, saddles and sources from the trivial basic set $\Omega_f$.  
We let $\Omega^{1u}$ be the set of saddle points $p\in \Omega^1$ for which there
is either a saddle point  $q\in(\Omega^{1}\setminus p)$ such that $W^u_{p}\cap W^s_{q}\neq\emptyset$ or a set $\Lambda\in\mathcal L^s$ such that $W^u_{p}\cap W^s_\Lambda\neq\emptyset$.
Next define  $\Omega^{1s}=\Omega^1\setminus\Omega^{1u}$. Note that the definitions of the sets $\Omega^{1s}$ and $\Omega^{1u}$ are not symmetric, but,  by the class $\Psi$ assumptions, if $p\in \Omega^{1u}$ then there exists no saddle point 
$q$ for which $W^s_{p}\cap W^u_{q}\neq\emptyset$ for some saddle point
$q$ and that there exists no set $\Lambda\in\mathcal L^u$ such that $W^s_{p}\cap W^u_\Lambda\neq\emptyset$.

\begin{figure}[ht]
\begin{center}\epsfig
{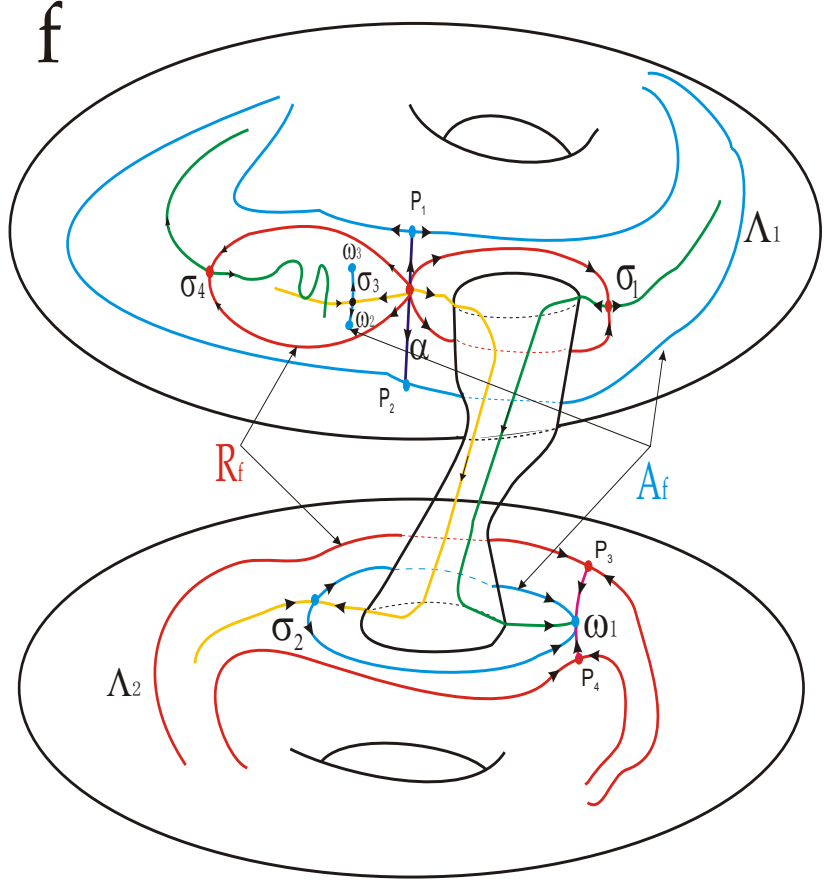} 
\caption{Attractor $A_f$ and  repeller $R_f$ for a diffeomorphism $f\in\Psi$}\label{AR}
\end{center}
\end{figure}

Let us set $$A_f=W^u_{\Omega^{1s}}\cup \Omega^0 \cup \mathcal L^s,~~~~~R_f=W^s_{\Omega^{1u}} \cup \Omega^2 \cup \mathcal L^u.$$
By construction the set $A_f$ is an attractor and $R_f$ is a repeller of $f$, see Figure \ref{AR}. 

\begin{figure}[ht]
\begin{center}\epsfig
{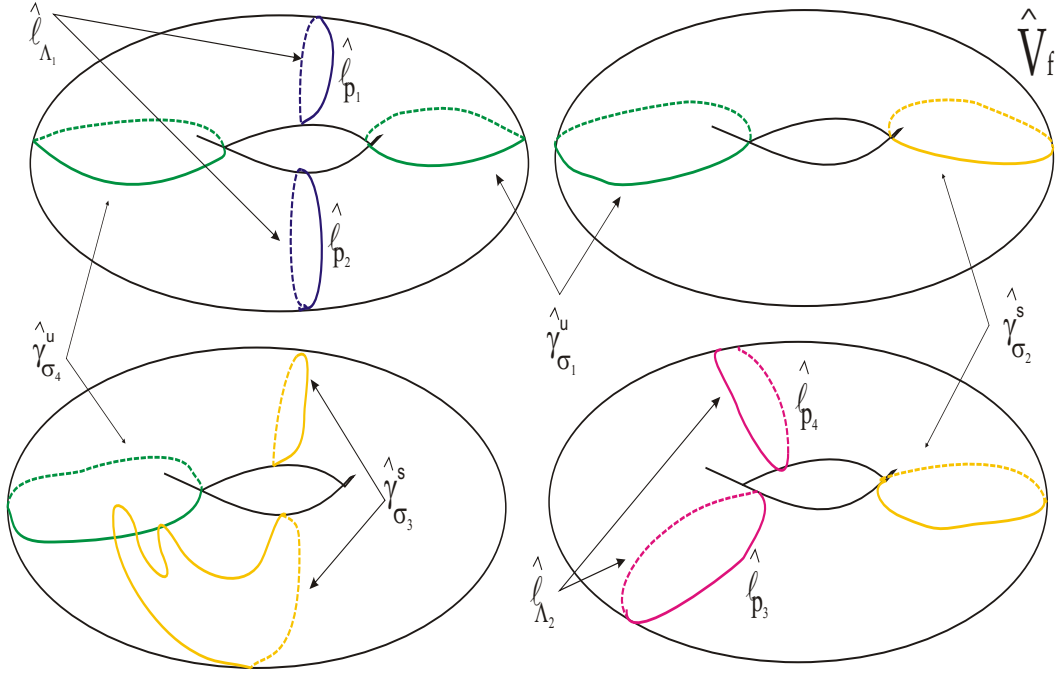} 
\caption{\small The orbit space $\hat V_f$ for the diffeomorphism $f$ from Figure \ref{AR} with the projections of the separatrices}\label{sch}
\end{center}
\end{figure}

Set $$V_f=M^2\setminus(A_f\cup R_f).$$ Let $k_{_f}$ be a minimal natural number such that each separatrix of saddle and boundary points is invariant with respect to $f^{k_{_f}}$. By construction the orbit space $\hat V_f=V_f/f^{k_{_f}}$ of the action of the diffeomorphism $f^{k_{_f}}$ on $V_f$ consists of a finite number of copies of the two-dimensional torus, and the natural projection $p_{_f}:V_f\to\hat V_f$ is a covering (see, for example \cite{MiPo2010}, Lemma 2.1). 
In Figure~\ref{sch} this construction is illustrated for the diffeomorphism shown in Figure~\ref{AR}. 
Set
$$\hat\phi_{_f}=p_{_{f}}fp_{_{f}}^{-1}:\hat{V}_{f}\to\hat{V}_{f}.$$ 

For each point $\sigma\in\Omega^{1s}$ set $\hat\gamma_{\sigma}=p_{_{f}}(W^{s}_{\sigma}\setminus\sigma)$ and for each point $\sigma\in\Omega^{1u}$ set $\hat\gamma_{\sigma}=p_{_{f}}(W^{u}_{\sigma}\setminus\sigma)$. By the construction $\hat\gamma_\sigma$ is a pair of circles.  Set  $$\hat\Gamma^s_f =\bigcup\limits_{\sigma\in\Omega^{1s}}\hat\gamma_{\sigma},~~\hat\Gamma^u_f =\bigcup\limits_{\sigma\in\Omega^{1u}}\hat\gamma_{\sigma}.$$
For each non-trivial basic set $\Lambda$ denote by $k_{_\Lambda}$ the number of periodic components of $\Lambda$ and by $g_{_\Lambda}$ the restriction of $f^{k_{_\Lambda}}$ on a periodic component of $N_\Lambda$. Set $\hat\ell_p=p_{_f}(\ell_p),~p\in P_\Lambda$,  $\hat\ell_\Lambda=\bigcup\limits_{p\in\Lambda}\hat\ell_p$ and $$\hat L^s_f=\bigcup\limits_{\Lambda\in\mathcal L^s}\hat\ell_\Lambda,~~\hat L^u_f=\bigcup\limits_{\Lambda\in\mathcal L^u}\hat\ell_\Lambda.$$ 

Set $\hat H_a=p_{_f}(H_a)$, $\lambda_{\hat H_a}=\lambda_a$ and $\mu_{\hat H_a}=\mu_a$ for $a\in\mathcal A$. Denote by $\hat H_f$ the union of all sets $\hat H_a$. For a connected component $\hat T$ of $\hat V_f$ let us set $\hat H_{\hat T}=\hat H_f\cap\hat T$. If the set $\hat H_{\hat T}$ is not empty let us choose a simple closed curve $\hat\beta_{\hat T}$ which intersects each curve from $\hat T\cap\hat\Gamma^s_f$ at exactly one point not being from $\hat{H}_f$ (such curve exists as $\hat T\cap\hat\Gamma^s_f$ is a set of disjoint non-contractible curves). 
Denote by $\beta_{\hat T}$ a connected component of the preimage $p_f^{-1}(\hat\beta_{\hat T})$ and by $K_{\hat T}$ an annulus on $V_f$ situated between $\beta_{\hat T}$ and $f^{k_f}(\beta_{\hat T})$. For an oriented path $\hat\nu\subset\hat T$ from a point $\hat x$ to a point $\hat y$ there is a unique lift $\nu\subset V_f$ with the start point $x= p^{-1}_f(\hat x)\cap K_{\hat T}$ (see, for example, \cite{Ko}). Then the end point of $\nu$ is situated in $f^{k_f\cdot k_{\hat\nu}}(K_{\hat T})$ for some integer $k_{\hat\nu}$.   

Let $\hat a_1,\hat a_2\in\hat H_{\hat T}$. If $\hat a_1$ and $\hat a_2$ belong to the same connected component of $\hat\Gamma^s_f$ then  denote by $\hat\nu_{\hat a_1,\hat a_2}$ a directed curve connecting the points $\hat a_1$ with $\hat a_2$ which is the part of curve from $\hat\Gamma^s_f$ oriented along the stable manifold. If $\hat a_1$ and $\hat a_2$ belong to different connected components $\hat\gamma^s_1,~\hat\gamma^s_2$ of  $\hat\Gamma^s_f$ then set $\hat z_1=\hat\gamma^s_1\cap\hat\beta_{\hat T},~\hat z_1=\hat\gamma^s_1\cap\hat\beta_{\hat T}$ and denote by $\hat\nu_{\hat a_1,\hat a_2}$ a directed curve connecting the points $\hat a_1$ with $\hat a_2$ consisting of a part of  $\hat\gamma^s_1$ oriented opposite the stable manifold, a part of curve $\hat\beta_{\hat T}$ connecting $\hat z_1$ with $\hat z_2$ and a part of  $\hat\gamma^s_2$ oriented along the stable manifold. 

For each point $b\in(H_f\cap K_{\hat T})$ let us calculate $\tau_b$ and set $\tau_{\hat b}=\tau_b$ for $\hat b=p_{_f}(b)$. For $\hat H_a$ is from $\hat H_f$ let us set $$\tau_{_{\hat H_a}}=\{\tau_{\hat b},~\hat b\in\hat H_a\}~~~and~~~\hat C_{_{\hat H_a}}=\{\lambda_{\hat H_a},~\mu_{\hat H_a},\tau_{_{\hat H_a}}\}.$$ Set 
$$\hat C_f=\{\hat C_{\hat H_{a}},\hat H_a\subset\hat H_f\}.$$ 

\begin{defi}[The scheme of a diffeomorphism] We call the set $S_f=(\hat{V}_{f}, \phi_{_f}, \hat\Gamma^s_f, \hat\Gamma^u_f,\hat{C}_{f}, \hat L^s_f, \hat L^u_f)$ a scheme of the diffeomorphism  $f\in\Psi$.
\end{defi}

\begin{defi}[Equivalence of schemes] The schemes 
$$S_f=(\hat{V}_{f}, \phi_{_f},\hat\Gamma^s_f, \hat\Gamma^u_f,\hat{C}_{f}, \hat L^s_f, \hat L^u_f)\mbox{ and }S_{f'}=(\hat{V}_{f'}, \phi_{_{f'}} \hat\Gamma^s_{f'}, \hat\Gamma^u_{f'},\hat{C}_{f'}, \hat L^s_{f'}, \hat L^u_{f'})$$ of diffeomorphisms  $f, f' \in \Psi$, respectively, are said to be equivalent if there exists an orientation-preserving homeomorphism $\hat\varphi:\hat V_f\to\hat V_{f'}$ such that :

1) $\hat\varphi\hat\phi_{_f}=\hat\phi_{_{f'}}\hat\varphi$;

2) $\hat\varphi(\hat\Gamma^s_f)=\hat\Gamma^s_{f'}$, $\hat\varphi(\hat\Gamma^u_f)=\hat\Gamma^u_{f'}$ and for each $\sigma\in\Omega^1$ there is a unique $\sigma'\in\Omega'^1$ such that  $\hat\varphi(\hat\gamma_{\sigma})=\hat\gamma_{\sigma'}$;

3) if $\hat H_{a'}=\hat\varphi(\hat H_a)$ then  $\frac{ln\left|\lambda_{\hat H_a}\right|}{ln\left|\mu_{\hat H_a}\right|}=\frac{ln\left|\lambda_{\hat H_{a'}}\right|}{ln\left|\mu_{\hat H_{a'}}\right|}$; 

4) if $\hat H_{a'}=\hat\varphi(\hat H_a)$ for $\hat H_a$ from $\hat H _f$ then 

$~~~4a)$ for any points $\hat a_1,~\hat a_2\in \hat H_a$ belonging to the  same connected component $\hat T$ of $\hat V_f$ we have $\left|\frac{\tau_{\hat a_2}}{\tau_{\hat a_1}}\right|^{\frac{1}{\ln\left|\mu_{\hat H_a}\right|}}=\left(\left|\frac{\lambda_{\hat H_{a'}}}{\mu_{\hat H_{a'}}}\right|^{k_{_{\hat\varphi(\hat\nu_{\hat a_1,\hat a_2})}}}\cdot\left|\frac{\tau_{\hat\varphi(\hat a_2)}}{\tau_{\hat\varphi(\hat a_1)}}\right|\right)^{\frac{1}{\ln\left|\mu_{\hat H_{a'}}\right|}}$, ;

$~~~4b)$ for any points $\hat a_1,~\hat a_2\in \hat H_a$ belonging to different connected components $\hat T_1,~\hat T_2$ of $\hat V_f$ there is a number $m_{_{\hat a_1,\hat a_2}}$ such that $\left(\frac{\tau_{\hat a_2}}{\tau_{\hat a_1}}\right)^{\frac{1}{\ln\left|\mu_{\hat H_a}\right|}}=\left(\left|\frac{\lambda_{\hat H_{a'}}}{\mu_{\hat H_{a'}}}\right|^{m_{_{\hat a_1,\hat a_2}}}\cdot\frac{\tau_{\hat\varphi(\hat a_2)}}{\tau_{\hat\varphi(\hat a_1)}}\right)^{\frac{1}{\ln\left|\mu_{\hat H_{a'}}\right|}}$;

5) if $H_a,~H_b$ are from $H_f$ and  $\hat a_1,~\hat a_2\in \hat H_a$, $\hat b_1,~\hat b_2\in \hat H_b$ such that $\hat a_1,\hat b_1\in\hat T_1,~\hat a_2,\hat b_2\in\hat T_2$ then the numbers $m_{_{\hat a_1,\hat a_2}}$, $m_{_{\hat b_1,\hat b_2}}$ satisfy the equality $m_{_{\hat b_1,\hat b_2}}=-k_{_{\hat\varphi(\hat\nu_{\hat a_1,\hat b_1})}}+m_{_{\hat a_1,\hat a_2}}+k_{_{\hat\varphi(\hat\nu_{\hat a_2,\hat b_2})}}$; 

6) $\hat\varphi(\hat L^s_f)=\hat L^s_{f'}$, $\hat\varphi(\hat L^u_f)=\hat L^u_{f'}$ and for each $\Lambda\in\mathcal L$ there is a unique $\Lambda'\in\mathcal L'$ such that  $\hat\varphi(\hat\ell_{_\Lambda})=\hat\ell_{_{\Lambda'}}$;

7) if $\hat\varphi(\hat\ell_{_\Lambda})=\hat\ell_{_{\Lambda'}}$ then there is an isomorphism $\psi_{_\Lambda}$ conjugating $T_{\bar {g}_{_\Lambda}}$ with $T_{\bar {g}^{\prime}_{\Lambda'}}$ for some $\bar {g}_{_\Lambda},\bar {g}^{\prime}_{\Lambda'}$ and such that $\hat\varphi(\hat\ell_{p})=\hat\ell_{_{\hat\psi_{\Lambda}(p)}}$. 
\label{eqv}
\end{defi}

\section{Separability of a 1-dimensional attractor (repeller) of a  diffeomorphism of a surface with a finite number moduli}
\label{classs}

{\bf Proof of Theorem \ref{1-sep}} 

We now prove that a 1-dimensional attractor of an A-diffeomorphism $f:M^2\to M^2$ with a finite number moduli is separable. 

\begin{demo} Let $\Lambda$ be an attractor  of an A-diffeomorphism $f:M^2\to M^2$ with a finite number of moduli. 
Let us prove that the three conditions of Definition \ref{ott} hold.

1) To prove item 1) of Definition \ref{ott} it suffices to prove that $W^u_{\Lambda'}\cap W^s_\Lambda=\emptyset$ holds for every non-trivial basic set $\Lambda'$ distinct from $\Lambda$. Suppose the contrary: there are points $x\in\Lambda,x'\in\Lambda^{\prime}$ such that $W^{s}_x\cap W^{u}_{x'}\neq\emptyset$. Since stable manifolds of the points of $\Lambda$ (unstable manifolds of the points of $\Lambda^{\prime}$) are $C^{1}$-close on compact sets, without loss of generality one can assume that the manifold $W^{s}_x$ contains no $s$-boundary periodic points of the basic set $\Lambda$ and the manifold $W^{u}_{x'}$ contains no $u$-boundary periodic points of the basic set $\Lambda'$. By Statement \ref{fimo}, the intersection $W^u_{\Lambda'}\cap W^s_\Lambda$ is transverse.   

Let $y\in(W^{s}_x\cap W^{u}_{x'})$. As $\Lambda$ and $\Lambda'$ have local structure of the product on interval by Cantor set then the point $y$ belongs to an adjacent interval $(a,b)^s\subset W^s_x$ which consists of the wandering points of the diffeomorphism $f$ and such that $a,b\in\Lambda$ and $W^u_a, W^u_b$ contain one $s$-boundary point each $p_a,p_b$ ($p_a=p_b$ if $W^u_a=W^u_b$). Denote by   $L_a$ ($L_b$) the connected component of the set $W^{u}_a\setminus a$ ($W^{u}_b\setminus b$) disjoint from the point  $p_a$ ($p_b$). Then the curve $l_{ab} = L_a\cup L_b\cup [a,b]^{s}$ bounds a domain $D_{ab}$. This domain is a continuous immersion of the open disk into the manifold $M^2$, all of its points are the wandering points of the diffeomorphism $f$ and the curve $l_{ab}$ is the boundary of $D_{ab}$ which is accessible from inside.

Denote by $W^{u*}_y$ the connected component of the set $W^{u}_y\setminus y$ disjoint from the point $x'$.  The transversality condition implies $W^{u*}_y\cap D_{ab} \neq \emptyset$.  On the other hand the component $W^{u*}_y$ contains a set which is dense in the periodic component of the set $\Lambda^{\prime}$. Therefore there are points in $W^{u*}_y$ disjoint from the domain $D_{ab}$. Then there is a point  $y'\in (a,b)^{s}$ distinct from the point $y$ and such that the arc $(y,y')^{u}\subset W^u_{x'}$ belongs to the domain  $D_{xy}$. Since for any point $\tilde a\in L_{p_a}$ there is a unique point  $\tilde  b\in L_{p_b}$  such that   $\tilde a\in W^{s}_{\tilde  x},\tilde x\in\Lambda$  and  $(\tilde  a,  \tilde b)^{s} \subset D_{ab}$ it follows that there is a point $\tilde x$ for which the arc $(\tilde  a,  \tilde b)^{s}$ is tangent to the arc $(y,y')^{u}$ and this contradicts the  transversality condition.

2) To prove the item 2) of Definition \ref{ott} it suffices to show that for every $s$-boundary point $p$ of the basic set $\Lambda$ there is no saddle point $\sigma$ from the trivial basic set of the diffeomorphism $f$ such that $W^u_\sigma\cap \ell_p\neq\emptyset$. If we assume the contrary then similarly to the proof of the item 1) we come to a  contradiction to the transversality condition.

3) Assuming the contrary in this case we come to a  contradiction to the transversality condition as well. 
\end{demo}  

\section{A proof of the classification theorem}

{\bf A proof of Theorem \ref{th1}}

It follows from the geometrical construction of the schemes and Statement \ref{mnogo} that diffeomorphisms $f,f'\in\Psi^*$ are topologically conjugate then  their schemes are equivalent. Let us show that if the schemes of diffeomorphisms $f,f'\in\Psi$ are equivalent, then the diffeomorphisms are topologically conjugate.

\begin{demo} Let $S_f=(\hat{V}_{f}, \phi_{_f},\hat\Gamma^s_f, \hat\Gamma^u_f,\hat{C}_{f}, \hat L^s_f, \hat L^u_f)$ and $S_{f'}=(\hat{V}_{f'}, \phi_{_{f'}} \hat\Gamma^s_{f'}, \hat\Gamma^u_{f'},\hat{C}_{f'}, \hat L^s_{f'}, \hat L^u_{f'})$ schemes of diffeomorphisms  $f, f' \in \Psi$, respectively, for which there exists an orientation-preserving homeomorphism $\hat\varphi:\hat V_f\to\hat V_{f'}$ with the properties 1)-7) of Definition \ref{eqv}. We divide the construction of a conjugating homeomorphism $h:M^2\to M^2$ such that $hf=f'h$ in to steps.

{\bf Step 1.} The existence of the homeomorphism $\varphi:\hat V_f\to\hat V_{f'}$ with property $\hat\varphi\hat\phi_{_f}=\hat\phi_{_{f'}}\hat\varphi$ implies that there exists a homeomorphism $\varphi: V_f\to V_{f'}$ that conjugates the restriction of the diffeomorphism $f$ to $V_f$ with the restriction of the diffeomorphism $f'$ to $V_{f'}$ and is such that $\hat\varphi=p_{_{f'}} \varphi p^{-1}_{_{f}}$. For each point $b\in(H_a\cap K_{\hat T})$ let us denote by $n(\hat b)$ an integer such that $\varphi(b)\in f'^{k_{f'}\cdot n(\hat b)}(K_{\hat\varphi(\hat T)})$. Due to condition 5) in Definition \ref{eqv}, we can suppose that $\varphi$ is chosen such that if $a_1\in(H_a\cap K_{\hat T_1})$ and $a_2\in(H_a\cap K_{\hat T_2})$ then $n(\hat a_2)-n(\hat a_1)=m_{_{\hat a_1,\hat a_2}}$. So we have a conjugating homeomorphism on the set $M^2\setminus(A_f\cup R_f)$.

Due to condition 2) in Definition \ref{eqv}, for any point $\sigma\in \Omega^{1\delta},\delta\in\{u,s\}$ there exists a point $\sigma'\in \Omega'^{1\delta}$ such that $\varphi(W^{\delta}_{\sigma}\setminus  \sigma)=W^{\delta}_{\sigma'}\setminus \sigma'$. Let us extend $\varphi$ on the set  $\Omega^1$ by setting $\varphi(\sigma)=\sigma'$. Due to condition 1) in Definition \ref{eqv}, $\varphi$ conjugates $f|_{W^{\delta}_{\Omega^{1\delta}}}$ with  $f'|_{W'^{\delta}_{\Omega'^{1\delta}}}$.

Due to condition 6) in Definition \ref{eqv}, for any basic set $\Lambda\in\mathcal L^\delta$ there exists a basic set $\Lambda'\in \mathcal L'^\delta$ such that $\varphi(\ell_{_\Lambda})=\ell_{_{\Lambda'}}$. Let us extend $\varphi$ by continuity on the set  $\bigcup\limits_{\Lambda\in\mathcal L}P_\Lambda$ of boundary points of non-trivial basic sets. Due to condition 1) in Definition \ref{eqv}, $\varphi$ conjugates $f|_{\bigcup\limits_{\Lambda\in\mathcal L}\ell_{_\Lambda}}$ with  $f'|_{\bigcup\limits_{\Lambda'\in\mathcal L'}\ell_{_{\Lambda'}}}$.

{\bf Step 2.} In this step we define homeomorphisms  $\varphi^{s}_{\Omega^{1u}}:W^{s}_{\Omega^{1u}} \to
W^{s}_{\Omega'^{1u}}$ and $\varphi^{u}_{\Omega^{1s}}:W^{u}_{\Omega^{1s}} \to
W^{u}_{\Omega'^{1s}}$ which conjugate $f|_{W^{s}_{\Omega^{1u}}}$ with $f'|_{W^{s}_{\Omega'^{1u}}}$ and $f|_{W^{u}_{\Omega^{1s}}}$ with $f'|_{W^{u}_{\Omega'^{1s}}}$, accordingly. 

Let $\sigma\in\Omega^{1u}$ and $\sigma'=\varphi(\sigma)$. Set $$\rho_\sigma^s=\frac{\ln|\lambda_{\sigma'}|}{\ln|\lambda_{\sigma}|}.$$ Let $\ell^s_\sigma$ be a stable  separatrix of $\sigma$. Denote by $\ell^s_{\sigma'}$ the stable  separatrix of $\sigma'=\varphi(\sigma)$ such that for a connected component $E$ of $U_\sigma\setminus W^u_\sigma$ containing $\ell^s_\sigma$ and a connected component $E'$ of $U_{\sigma'}\setminus W^u_{\sigma'}$ containing $\ell^s_{\sigma'}$ the following condition holds $\varphi(E\setminus W^u_\sigma)\cap E'\neq\emptyset$. Let us define a homeomorphism $\varphi_{\ell^s_\sigma}:\ell^{s}_{\sigma} \to \ell^{s}_{\sigma'}$ by the following way. For a point $t\in \ell^{s}_{\sigma}$ such that $t^u=(0,t^{u}_{y})$ let us set  $\varphi_{\ell^s_\sigma}(t)=t'$ where ${t'}^u=(0,{t'}^{u}_{y})$, $$|{t'}^{u}_{y}|=|t^{u}_{y}|^{\rho_\sigma^s}.$$ It is easy to verify that $\varphi_{\ell^s_\sigma}$ conjugates the diffeomorphisms $f^{k_\sigma}\vert_{\ell^{s}_{\sigma}}$ and $f^{\prime k_\sigma}\vert_{\ell^{s}_{\sigma'}}$. 

Due to property 2) of Definition \ref{eqv} we get $k_\sigma=k_{\sigma'}$. Then for each $k=0,\dots,k_\sigma$ we can define a homeomorphism $\varphi_{\ell^s_{f^k(\sigma)}}:\ell^s_{f^k(\sigma)}\to \ell^s_{f^{\prime k}(\sigma')}$ by the formula $\varphi_{\ell^s_{f^k(\sigma)}}(x)=f^{\prime k}(\varphi_{\ell^s_\sigma}(f^{-k}(x)))$ for each $x\in \ell^s_{f^k(\sigma)}$. Doing a similar construction for all saddle periodic orbits of the set $\Omega^{1u}$ we get the sought conjugating homeomorphism $\varphi^{s}_{\Omega^{1u}}$.

Now let $\sigma\in \Omega^{1s}$ and $\sigma'=\varphi(\sigma)$. Set $$\rho_\sigma^u=\frac{\ln|\mu_{\sigma'}|}{\ln|\mu_{\sigma}|}.$$
Similar to the construction above for corresponding separatrices $\ell^u_\sigma,~\ell^u_{\sigma'}$ we define a homeomorphism $\varphi_{\ell^u_\sigma}:\ell^{u}_{\sigma} \to \ell^{u}_{\sigma'}$ by the following way. For a point $t \in \ell^{u}_{\sigma}$ such that $t^s=(t^{s}_{x},0)$ let us set  $\varphi_{\ell^u_\sigma}(t)=t'$ where ${t'}^s=({t'}^{s}_{x},0)$, $$|{t'}^{s}_{x}|=c_{\ell^u_\sigma}|t^{s}_{x}|^{\rho_\sigma^u},$$ where $c_{\ell^u_\sigma}=\frac{\left|\frac{\lambda_{\hat H_{a'}}}{\mu_{\hat H_{a'}}}\right|^{n(\hat a)}\cdot\vert \tau_{\hat a'}\vert}{\vert\tau_{\hat a}\vert^{\rho_{\sigma}^u}}$ if there is a point $\hat a\in\hat\ell^u_\sigma\cap\hat H_f$ and equals 1 in the opposite case. As above it is possible to verify that $\varphi^{u}_{\sigma}$ conjugates the diffeomorphisms $f^{k_\sigma}\vert_{\ell^{u}_{\sigma}}$ and $f^{\prime k_\sigma}\vert_{\ell^{u}_{\sigma'}}$. For each $k=0,\dots,k_\sigma$ we can define a homeomorphism $\varphi_{\ell^u_{f^k(\sigma)}}:\ell^u_{f^k(\sigma)}\to \ell^u_{f^{\prime k}(\sigma')}$ by the formula $\varphi_{\ell^u_{f^k(\sigma)}}(x)=f^{\prime k}(\varphi_{\ell^u_\sigma}(f^{-k}(x)))$ for each $x\in \ell^u_{f^k(\sigma)}$. 
Doing a similar construction for all saddle periodic orbits of the set $\Omega^{1s}$ we get the sought conjugating homeomorphism $\varphi^{u}_{\Omega^{1s}}$.

{\bf Step 3.} In this step we construct a homeomorphism $\varphi_{_{\mathcal L^{s}}}:\mathcal L^{s}\to{\mathcal L}^{\prime s}$ ($\varphi_{_{\mathcal L^{u}}}:\mathcal L^{u}\to{\mathcal L}^{\prime u}$) which conjugates $f|_{\mathcal L^s}$ with $f'|_{{\mathcal L}^{\prime s}}$ ($f|_{\mathcal L^u}$ with $f'|_{\mathcal L^{\prime u}}$). Let us construct $\varphi_{_{\mathcal L^{s}}}$, the construction of $\varphi_{_{\mathcal L^{u}}}$ is similar.

Let $\Lambda$ be a one-dimensional attractor of $f$ and $L$ be one from $k_{_\Lambda}$ periodic components of $\Lambda$. Then $L$ is a one-dimensional attractor of the diffeomorphism $g=f^{k_{_\Lambda}}$ with the unique periodic component. Denote by $L'$ a periodic component of $\Lambda'$ such that $\varphi$ sends the boundary points of $L$ to the boundary points of $L'$. Then $L'$ is a one-dimensional attractor of the diffeomorphism $g'=f'^{k_{_{\Lambda'}}}$ with the unique periodic component. Due to conditions 1) and 6) in Definition \ref{eqv}, $k_{_\Lambda}=k_{_{\Lambda'}}$. Due to condition 7) in Definition \ref{eqv}, there is an isomorphism $\psi_{_\Lambda}$ conjugating $T_{\bar g_{_L}}$ with $T_{\bar g'_{L'}}$ for some $\bar g_{_L},\bar g'_{L'}$ and such that $\varphi(p)=\hat\psi_{_\Lambda}(p)$ for each point $p$ from the set $P_L$ of the boundary points of $L$. Statement \ref{atr.l8} implies that there is a homeomorphism $\bar{\varphi}_{L}:\bar L\to\bar L'$ such that $\bar{g}'_{L'} \bar{\varphi}_{L} \vert_{\bar{L}}=\bar{\varphi}_{L} \bar{g}_{_L} \vert_{\bar{L}}$. Set
$$\varphi_L=p_{_{N_{L'}}}\bar{\varphi}_{L}p^{-1}_{N_L}:L\to L^\prime.$$ 
Then
$\varphi_L g_L\vert_L=g'_{L^{\prime}}\varphi_L\vert_L$ and $\varphi_L(p)=\varphi(p)$ for each $p\in P_L$. Define $\varphi_\Lambda:\Lambda\to\Lambda'$ by the  formula $\varphi_\Lambda(v)=f^{\prime k}(\varphi_L(f^{-k}(v)))$ where $f^k(v)\in L$ for $v\in\Lambda$. Doing a similar construction for all attractors of the set $\mathcal L^s$ we get the sought conjugating homeomorphism $\varphi_{_{\mathcal L^{s}}}$.

{\bf Step 4.} In this step we modify the homeomorphism $\varphi|_{W^s_{\mathcal L^s}\setminus\mathcal L^s}~(\varphi|_{W^u_{\mathcal L^u}\setminus\mathcal L^u})$ up to $\tilde h_{\mathcal L^s}:W^s_{\mathcal L^s}\setminus\mathcal L^s\to W^s_{\mathcal L^{\prime s}}\setminus\mathcal L^{\prime s}~(\tilde h_{\mathcal L^u}:W^s_{\mathcal L^u}\setminus\mathcal L^u\to W^s_{\mathcal L^{\prime u}}\setminus\mathcal L^{\prime u})$ which extended to $\mathcal L^s~(\mathcal L^u)$ by $\varphi_{_{\mathcal L^{s}}}~(\varphi_{_{\mathcal L^{u}}})$. Let us construct $\tilde h_{_{\mathcal L^{s}}}$, the construction of $\tilde h_{_{\mathcal L^{u}}}$ is similar.

Let $\Lambda$ be a one-dimensional attractor of $f$. Let us modify $\varphi|_{W^s_\Lambda\setminus \Lambda}$ up to $\tilde h_\Lambda:W^s_\Lambda\setminus\Lambda\to W^s_{\Lambda'}\setminus\Lambda'$ which extended to $\Lambda$ by $\varphi_{_\Lambda}$ and conjugates $f|_{W^s_\Lambda\setminus \Lambda}$ with $f'|_{W^s_{\Lambda'}\setminus \Lambda'}$.  

Let $L$ be a periodic component of $\Lambda$ as in Step above. Denote by $B_L$ the set of all bunches of $L$ and will use further the denotations of section \ref{Lam}.  Let us consider the closed curve $L_b$ of the bunch $b\in B_L$ and enumerate the separatrices $l_1,\dots,l_m$ of all saddle points and boundary points which intersect $L_b$ due to some orientation on $L_b$. Set $b'=\varphi_L(b)$. Due to items 2) and 6) of Definition \ref{eqv} we have that the separatrices $\varphi(l_1),\dots,\varphi(l_m)$ intersect $L_{b'}$ in order. If $\varphi_L|_{P_L}=\varphi|_{P_L}$ then for each $j\in\{1,\dots,r_b\}$ there is a homeomorphism $h^s_j:[x_{2j}, x_{2j+1}]^{s}\to[\varphi_L(x_{2j}), \varphi_L(x_{2j+1})]^{s}$ such that $h^s_j([x_{2j}, x_{2j+1}]^{s}\cap l_\mu)=[\varphi_L(x_{2j}), \varphi_L(x_{2j+1})]^{s}\cap\varphi(l_\mu)$ for each $\mu\in\{1,\dots,m\}$ and $h^s_j(x_{2j})=\varphi_L(x_{2j})$. Set $I^s_b=\bigcup\limits_{j=1}^{r_b}[x_{2j}, x_{2j+1}]^{s}$, $I^s_{b'}=\bigcup\limits_{j=1}^{r_{b'}}[\varphi_L(x_{2j}),\varphi_L(x_{2j+1})]^{s}$ and denote by $h^s_b:I^s_b\to I^s_{b'}$ a homeomorphism which composed from $h^s_j,j\in\{1,\dots,r_b\}$. Set $I^s_L=\bigcup\limits_{b\in B_L}I^s_b$, $I^s_{L'}=\bigcup\limits_{b'\in B_{L'}}I^s_{b'}$ and denote by $h^s_L:I^s_L\to I^s_{L'}$ a homeomorphism which composed from $h^s_b,b\in B_L$.    

\begin{figure}[ht]
\begin{center}\includegraphics[width=10cm]{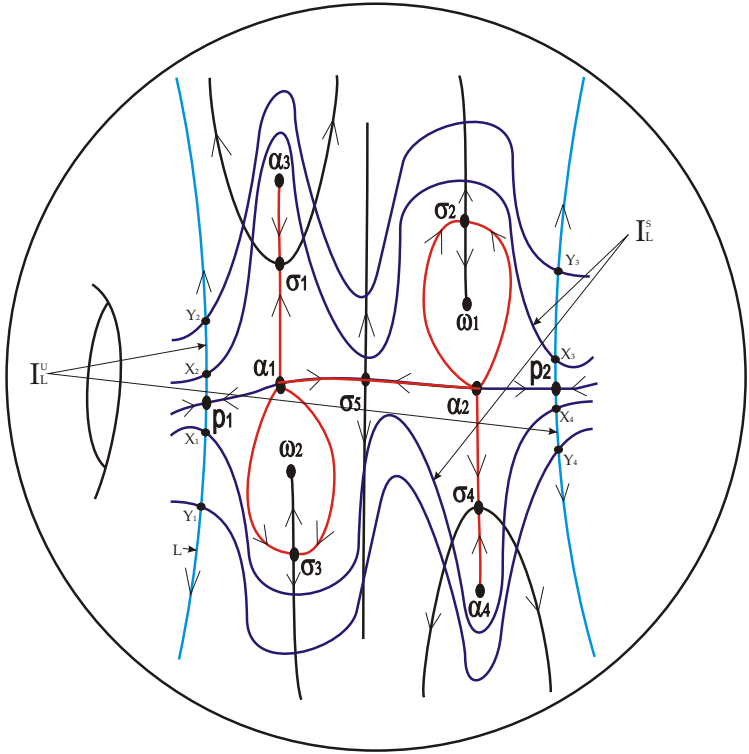}\end{center}
\caption{Illustration to the Step 4}\label{22b}
\end{figure}   

Denote by $y_{2j-1},y_{2j}\in W^u_{p_j}$ the intersection points of $f(I^s_b)$ with $W^u_{p_j}$ such that $p_j\notin[x_{2j},y_{2j}]^u$ (see figure \ref{22b}). Set $I^u_b=\bigcup\limits_{i=j}^{r_b}[x_{2j}, y_{2j}]^{u}$, $I^u_L=\bigcup\limits_{b\in B_L}I^u_b$, and $I^u_{L'}=h_L(I^u_L)$. Set $h^u_L=\varphi_L|_{I^u_L}:I^u_L\to I^u_{L'}$. Let $\Pi_L~(\Pi_{L'})$ be the closure of the set $W^s_{I^u_L}\setminus L~(W^s_{I^u_{L'}}\setminus L')$. Let us construct on $\Pi_L$ a pair of transverse one-dimensional foliation $F^s_L,~F^u_L$ with the following properties:

a) each leaf of $F^s_L$ is a connected component of the intersection stable manifold of a point from $L$ with $\Pi_L$;

b) each leaf of $F^u_L$ is a segment $[x,y]$ with the boundary points $x,~y$ such that $x\in[x_{2j},x_{2j+1}]^s,~y\in[y_{2j},y_{2j+1}]^s$;

c) if $[x,y]$ belongs to $F^u_L$ then $[g^{-1}(y),g(x)]$ also belongs to ot;

d) each connected component of intersection of the separatrices of the saddle points with $\Pi_L$ is a leaf of $F^u_L$. 

For each point $z\in I^u_L$ denote by $F^s_{L,z}$ a leaf of the foliation $F^s_L$ passing through the point $z$. For each point $x\in I^s_L$ denote by $F^u_{L,x}$ a leaf of the foliation $F^u_L$ passing through the point $x$. Let us construct similar foliations $F^s_{L'},~F^u_{L'}$ on $\Pi_{L'}$ and define a homeomorphism $h_{\Pi_L}:\Pi_L\to\Pi_{L'}$ by the formula $$h_{\Pi_L}(F^s_{L,z}\cap F^u_{L,x})=F^s_{L,h^u_L(z)}\cap F^u_{L,h^s_L(x)}.$$
Notice that $\Pi_L\setminus B_L$ is a fundamental domain of $f$ restriction on $W^s_\Lambda\setminus \Lambda$. Then for each point $w\in(W^s_\Lambda\setminus \Lambda)$ there is $k\in\mathbb Z$ such that $f^{k}(w)\in\Pi_L$. As $h_{\Pi_L}$ conjugates $g|_{\Pi_L}$ with $g'|_{\Pi_{L'}}$ then we can extend $h_{\Pi_L}$ up to $\tilde h_{\Lambda}:W^s_\Lambda\setminus \Lambda\to W^s_{\Lambda'}\setminus \Lambda'$ conjugating $f$ and $f^{\prime}$ by the formula $$\tilde h_{\Lambda}(w)=f^{\prime k}(h_{\Pi_L}(f^{-k}(w))).$$ 

Doing a similar construction for all attractors of the set $\mathcal L^s$ we get the sought conjugating homeomorphism $\tilde h_{_{\mathcal L^{s}}}$.

Denote by $\varphi_1:V_f\to V_{f'}$ a homeomorphism given by the formula
$${\varphi}_{1}(z)=
\begin{cases}\tilde h_{\mathcal L^{s}}(z), z\in (W^s_{\mathcal L^s}\setminus\mathcal L^s); \cr \tilde h_{\mathcal L^{u}}(z), z\in (W^u_{\mathcal L^u}\setminus\mathcal L^u); \cr
{\varphi}(z),z\in V_f\setminus(W^s_{\mathcal L^s}\cap W^u_{\mathcal L^u}).\cr\end{cases}$$
Set $\hat\varphi_1=p_{_{f'}}\varphi_1p_f^{-1}:\hat V_f\to\hat V_{f'}$.

{\bf Step 5.} In this step we modify the homeomorphism $\varphi_1|_{W^s_{\mathcal L^s}\setminus\mathcal L^s}~(\varphi_1|_{W^u_{\mathcal L^u}\setminus\mathcal L^u})$ up to $h_{\mathcal L^s}:W^s_{\mathcal L^s}\setminus\mathcal L^s\to W^s_{\mathcal L^{\prime s}}\setminus\mathcal L^{\prime s}~(h_{\mathcal L^u}:W^s_{\mathcal L^u}\setminus\mathcal L^u\to W^s_{\mathcal L^{\prime u}}\setminus\mathcal L^{\prime u})$ which extended to $\mathcal L^s~(\mathcal L^u)$ by $\varphi_{_{\mathcal L^{s}}}~(\varphi_{_{\mathcal L^{u}}})$ and extended by $\varphi^s_{\Omega^{1u}}~(\varphi^u_{\Omega^{1s}})$ on $cl~(W^s_{\mathcal L^s})\setminus(\mathcal L^s\cup\Omega^2)~(cl~(W^u_{\mathcal L^u})\setminus(\mathcal L^u\cup\Omega^0))$. Let us construct $h_{_{\mathcal L^{s}}}$, the construction of $h_{_{\mathcal L^{u}}}$ is similar. 

By Theorem \ref{1-sep}, each non-trivial attractor of the diffeomorphism $f$ is separable. Then there is a set $\Sigma^u\subset\Omega^{1u}$ such that $cl~(W^s_{\mathcal L^s})\setminus(\mathcal L^s\cup\Omega^2)=W^s_{\Sigma^u}$. Let $\sigma\in\Sigma^u$. Set $h^s_\sigma=\varphi^s_{\Omega^{1u}}|_{W^s_\sigma}:W^s_\sigma\to W^s_{\sigma'}$ and $h^u_\sigma=\varphi_1|_{W^u_\sigma}:W^u_\sigma\to W^u_{\sigma'}$. In an $f^{k_\sigma}$-invariant neighbourhood $N_\sigma$ of $\sigma$ let us construct a pair of transverse $f^{k_\sigma}$-invariant foliations $G^s_\sigma,G^u_\sigma$ with the following properties:

a) $W^s_\sigma\in G^s_\sigma$, $W^u_\sigma\in G^u_\sigma$;

b) if $W^u_\sigma\cap W^s_L\neq\emptyset$ for some periodic component $L$ of a non-trivial attractor then each connected component of $W^s_x\cap N_\sigma,x\in L$ is a leaf of $G^s_\sigma$.     

For each point $z_u\in W^u_\sigma$ denote by $G^s_{\sigma,z_u}$ a leaf of the foliation $G^s_\sigma$ passing through the point $z_u$. For each point $z_s\in W^s_\sigma$ denote by
$G^u_{\sigma,z_s}$ a leaf of the foliation $G^u_\sigma$ passing through the point $z_s$.
Let us construct similar foliations $G^s_{\sigma'},~G^u_{\sigma'}$ on $N_{\sigma'}$ and define a homeomorphism $h_{N_\sigma}:N_\sigma\to N_{\sigma'}$ by the formula $$h_{N_\sigma}(G^s_{\sigma,z_u}\cap G^u_{\sigma,z_s})=G^s_{\sigma',h^u_\sigma(z_u)}\cap G^u_{\sigma',h^s_\sigma(z_s)}.$$
Then in some tubular neighbourhood $N(\hat\gamma_\sigma)$ of $\hat\gamma_\sigma$ a map $\hat h_{N_\sigma}$ is well-defined by the formula $\hat h_{N_\sigma}=p_{_{f'}}h_{N_\sigma}p_{_f}^{-1}$. Chose a tubular neighbourhood $\tilde N(\hat\gamma_\sigma)$ of $\hat\gamma_\sigma$ such that $N(\hat\gamma_\sigma)\subset \tilde N(\hat\gamma_\sigma)$, $\hat h_{N_\sigma}(N(\hat\gamma_\sigma))\subset\hat\varphi_1( \tilde N(\hat\gamma_\sigma))$ and the set $Q=cl~(\tilde N(\hat\gamma_\sigma)\setminus  N(\hat\gamma_\sigma))$, ${Q}'=cl~(\hat\varphi_1( \tilde N(\hat\gamma_\sigma))\setminus\hat h_{N_\sigma}(N(\hat\gamma_\sigma)))$ are two-dimensional annulus.  Then there is a homeomorphism  $\hat\varphi_{\hat Q}:\hat Q\to \hat Q'$ such that  $\hat\varphi_{\hat Q}\vert_{\partial N(\hat\gamma_\sigma)}=\hat h_{N_\sigma}$ and $\hat\varphi_{\hat Q}\vert_{\partial \tilde N(\hat\gamma_\sigma)}=\hat\varphi_{1}$. As the homeomorphisms $\varphi_1$ and $h_{N_\sigma}$ send leaves of  the foliation $W^s_x,x\in\mathcal L^s$ to leaves of  the foliation $W^s_{x'},x'\in\mathcal L^{\prime s}$ and are coincide on $W^u_\sigma$ then we can construct $\hat\varphi_{\hat Q}$ such that its lift sends leaves of  the foliation $W^s_x,x\in\mathcal L^s$ to leaves of  the foliation $W^s_{x'},x'\in\mathcal L^{\prime s}$.  

Denote by $\hat\varphi_{\hat\gamma_\sigma}:\hat V_f\to\hat V_{f'}$ a homeomorphism given by the formula
$$\hat{\varphi}_{\hat\gamma_\sigma}(\hat z)=
\begin{cases} \hat\phi^k_{f'}(\hat h_{N_\sigma}(\phi^{-k}_f(\hat z))), \hat z\in\hat\phi^k_f(N(\hat\gamma_\sigma)); \cr \hat\phi^k_{f'}(\hat\varphi_{\hat Q}(\phi^{-k}_f(\hat z))), \hat z\in\hat\phi^k_f(\hat Q); \cr
\hat\varphi_1(\hat z),\hat z\in (\hat V_f\setminus\tilde N(\hat\gamma_\sigma)).\cr\end{cases}$$
Denote by $\varphi_{\hat\gamma_\sigma}$ a lift of $\hat\varphi_{\hat\gamma_\sigma}$ coinciding with $\varphi_1$ on $V_f\setminus p_f^{-1}(\tilde N(\hat\gamma_\sigma))$.
Doing in series a similar construction for all saddle periodic orbits of the set $\Sigma^{u}$ we get a homeomorphism $\varphi_{\Sigma^{u}}:V_f\to V_{f'}$. Also we construct a homeomorphism $\varphi_{\Sigma^{s}}:V_f\to V_{f'}$.

Denote by $\varphi_2:V_f\to V_{f'}$ a homeomorphism given by the formula
$${\varphi}_{2}(z)=
\begin{cases} \varphi_{\Sigma^u}(z), z\in (W^s_{\mathcal L^s}\setminus\mathcal L^s); \cr \varphi_{\Sigma^s}(z), z\in (W^u_{\mathcal L^u}\setminus\mathcal L^u); \cr
\hat\varphi_1(z),z\in V_f\setminus(W^s_{\mathcal L^s}\cap W^u_{\mathcal L^u}).\cr\end{cases}$$
Set $\hat\varphi_2=p_{_{f'}}\varphi_2p_f^{-1}:\hat V_f\to\hat V_{f'}$.

{\bf Step 6.}  Let $\sigma\in\Omega^1$. For a point $x \in {U}_{\sigma}$ denote by $\mathcal{F}^u_{\sigma,x}$ ($\mathcal{F}^s_{\sigma,x}$) the unique leaf of  $\mathcal{F}^u_{\sigma}$ ($\mathcal{F}^s_{\sigma}$) that passes through the point $x$. Define projections ${\pi}^{u}_{\sigma}:{U}_{\sigma}\to W^{s}_{\sigma}$
(${\pi}^{s}_{\sigma}: {U}_{\sigma} \to W^{u}_{\sigma}$) along the leaves of the foliation  $\mathcal{F}^u_{\sigma}$ ($\mathcal{F}^s_{\sigma}$) as follows:
${\pi}^{u}_{\sigma}(x)=\mathcal{F}^u_{\sigma,x}\cap W^{s}_{\sigma}$ (${\pi}^{s}_{\sigma}(x)=
\mathcal{F}^s_{\sigma,x}\cap W^{u}_{\sigma}$).
 
Let $a \in W^{s}_{\sigma^s_a} \cap W^{u}_{\sigma^u_a}$ be a point of one-sided tangency and $a'=\varphi(a)$. Set $l_a=\psi^{-1}_{\sigma^u_a}(\{(x,y) \in U_{\mu_{\sigma^u_a},\lambda_{\sigma^u_a}}: x=a^u_x\}) \cap U_{a}$, $l_{a'}=\psi^{-1}_{\sigma^u_{a'}}(\{(x,y) \in U_{\mu_{\sigma^u_{a'}},\lambda_{\sigma^u_{a'}}}: x={a'}^u_x\}) \cap U_{a'}$. 
Set $L_\mathcal A=\bigcup\limits_{a\in\mathcal A}l_a$ and $L_{\mathcal A'}=\bigcup\limits_{a'\in\mathcal A'}l_{a'}$. In this step we construct a  
homeomorphism $\varphi_{_{L_\mathcal A}}:L_\mathcal A\to L_{\mathcal A'}$ which conjugates  $f|_{L_\mathcal A}$ with $f'|_{L_{\mathcal A'}}$ and extended to $W^s_{\Omega^{1u}}$ by $\varphi^s_{\Omega^{1u}}$ and to $W^u_{\Omega^{1s}}$ by $\varphi^u_{\Omega^{1s}}$. 

Define a homeomorphism $\varphi_{l_a}:l_a\to l_{a'}$ by the formula $\varphi_{l_a}(z)=z'= ((\pi^u_{\sigma^u_{a'}})^{-1}(\varphi^{s}_{\sigma^u}(\pi^u_{\sigma^u_a}(z)))) \cap l_{a'}$. Set $L_{a}=\bigcup\limits_{n\in\mathbb{Z}}f^{kn}(l_a)$ and  $L_{a'}=\bigcup\limits_{n\in\mathbb{Z}}f'^{kn}(l_{a'})$, where $k$ is the period of unstable separatrix containing $a$. Define a homeomorphism $\varphi_{L_{a}}:L_{a} \to L_{a'}$ by the formula $\varphi_{L_{a}}(z)=z'=f^{\prime kn}(\varphi_{l_a}(f^{-kn}(z)))$ for each point $z\in f^{kn}(l_a)$. Set  $\mathcal{E}_{a}=W^{u}_{\sigma^s_a} \cup W^{s}_{\sigma^u_a} \cup L_a$ and  $\mathcal{E}_{a'}=W^{u}_{\sigma^s_{a'}} \cup W^{s}_{\sigma^u_{a'}} \cup L_{a'}$. Denote by  $\varphi_{a}:\mathcal{E}_{a}\to\mathcal{E}_{a'}$ a map, coinciding with  $\varphi^{u}_{\sigma^s_a}$ on $W^{u}_{\sigma^s_a}$, with $\varphi^{s}_{\sigma^u_a}$ on $W^{s}_{\sigma^u_a}$ and with $\varphi_{L_{a}}$ on $L_a$. Using condition  2) of Definition \ref{eqv} it is possible to verify that $\varphi_{a}$ is a homeomorphism (see \cite{MiPo2010} for details).

Denote by $A\subset\mathcal{A}$ a set of such points that any two from their are not belonging to the same orbit of $f$ and $\bigcup\limits_{n\in\mathbb Z}f^n(A)=\mathcal{A}$. Set $A'=\varphi(A)$, $\mathcal{E}_{\mathcal{A}}= \bigcup \limits_{a \in {A}} \mathcal{E}_{a}$ and $\mathcal{E}_{\mathcal{A}'}= \bigcup \limits_{a' \in {A}'} \mathcal{E}_{a'}$. Let us define a map $\varphi_{\mathcal{A}}:\mathcal{E}_{\mathcal{A}} \to \mathcal{E}_{\mathcal{A}'}$ as coinciding with $\varphi_{a}$ on each set $\mathcal{E}_{a}$. 
Using condition 4) of Definition \ref{eqv} it is possible to verify that  $\varphi_{\mathcal{A}}$ is a homeomorphism (see \cite{MiPo2010} for details). 

{\bf Step 7.} In the neighborhood $U_a$ of a point $a\in A$ define foliations  $\mathcal{F}^{u}_{a}$ and $\mathcal{F}^{s}_{a}$ by the following way. The leaves of $\mathcal{F}^{u}_{a}$ are coincide with the leaves of $\mathcal{F}^{u}_{\sigma^u_a}\cap U_a$. In the neighborhood $\psi_{\sigma^u_a}(U_a)$ the curve $\psi_{\sigma^u_a}(W^s_{\sigma^s_a})$ has the equation $q(x)=Q(x-a^u_{x})^n + o((x-a^u_{x})^n)$, where $\frac{o((x-a^u_{x})^n)}{(x-a^u_{x})^n} \to 0$ for $x \to a^u_{x}$. Set $\mathcal{F}^s_{a}=\psi^{-1}_{\sigma^u_a}(\bigcup\limits_{c\in\mathbb
R}\{(x,y)\in U_{\mu_{\sigma}, \lambda_{\sigma}}~:~y=q(x)+c\})\cap U_a$. Hence, in a neighborhood $U_a$ the leaves of ${\mathcal{F}}^{u}_{a}$ are transverse to the leaves of ${\mathcal{F}}^{s}_{a}$ on the set $U_a\setminus l_a$ and have tangency along the curve  $l_a$. Set $U_{A}=\bigcup\limits_{a\in A}U_a$, $U_{\mathcal A}=\bigcup\limits_{n\in\mathbb Z}f^n(U_A)$, $\mathcal{F}^{u}_{A}=\bigcup\limits_{a\in A}{\mathcal{F}}^{u}_{a}$, $\mathcal{F}^{u}_{\mathcal A}=\bigcup\limits_{n\in\mathbb Z}f^n({\mathcal{F}}^{u}_{A})$, $\mathcal{F}^{s}_{A}=\bigcup\limits_{a\in A}{\mathcal{F}}^{s}_{a}$ and  $\mathcal{F}^{s}_{\mathcal A}=\bigcup\limits_{n\in\mathbb Z}f^n({\mathcal{F}}^{s}_{A})$. The similar foliations $\mathcal{F}^{u}_{\mathcal A'}$ and $\mathcal{F}^{s}_{\mathcal A'}$ let us construct in the neighbourhood $U_{\mathcal A'}$ of the set $\mathcal {A'}$. 

Let $d$ be a point of the heteroclinic intersection of the manifolds $W^s_{\sigma_d^s}\cap W^u_{\sigma_d^u}$ be not belonging to the set $\mathcal A$. Denote $U_d$ a connected component of the set $U_{\sigma_d^s}\cap U_{\sigma_d^u}$ which contains $d$. Define foliations $\mathcal{F}^{u}_{d}$ a $\mathcal{F}^{s}_{d}$ by the following way: $\mathcal{F}^{u}_{d}=\mathcal{F}^{u}_{\sigma_d^u}\cap U_d$ and $\mathcal{F}^{s}_{d}=\mathcal{F}^{s}_{\sigma_d^s}\cap U_d$. Let  $\mathcal D$ ($\mathcal D'$) be the set of all heteroclinic points of $f$ ($f'$)  not belonging to $\mathcal A$, $D\subset\mathcal{D}$  the set of points such that any two of them do not belong to the same orbit of the diffeomorphism $f$ and $\bigcup\limits_{n\in\mathbb Z}f^n(D)=\mathcal{D}$. Set   $U_{D}=\bigcup\limits_{d\in D}U_d$, $U_{\mathcal D}=\bigcup\limits_{n\in\mathbb Z}f^n(U_D)$, $\mathcal{F}^{u}_{D}=\bigcup\limits_{d\in D}{\mathcal{F}}^{u}_{d}$, $\mathcal{F}^{u}_{\mathcal D}=\bigcup\limits_{n\in\mathbb Z}f^n({\mathcal{F}}^{u}_{D})$, $\mathcal{F}^{s}_{D}=\bigcup\limits_{d\in D}{\mathcal{F}}^{s}_{d}$ and  $\mathcal{F}^{s}_{\mathcal D}=\bigcup\limits_{n\in\mathbb Z}f^n({\mathcal{F}}^{s}_{D})$.  Construct the similar foliation $\mathcal{F}^{u}_{\mathcal D'}$ and $\mathcal{F}^{s}_{\mathcal D'}$ in the neighbourhood $U_{\mathcal D'}$ of the set  $D'=\varphi(D)$.

\begin{figure} \begin{center} \epsfig
{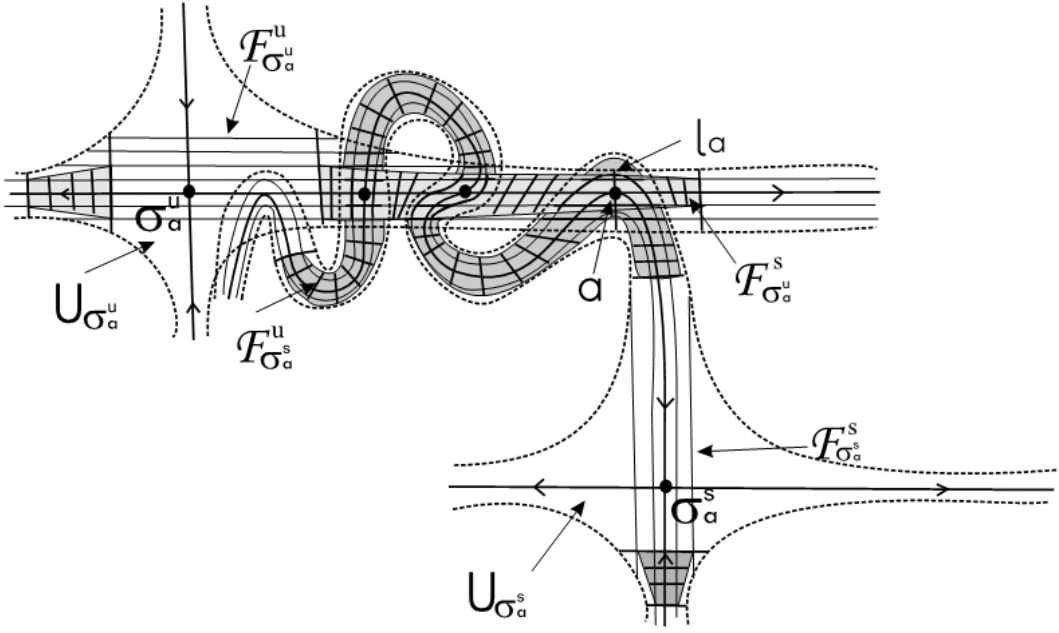} \caption{Construction of foliations} \label{conj} \end{center} 
\end{figure}

Set $L_{\mathcal A}=\bigcup\limits_{a\in A}L_a$. For a point $\sigma\in\Omega_f$ define foliations $\tilde{\mathcal{F}}^s_{\sigma}$ and $\tilde{\mathcal{F}}^u_{\sigma}$ transverse to each  other everywhere except $L_{\mathcal A}$ by the following way. The foliation  $\tilde{\mathcal{F}}^s_{\sigma}$ ($\tilde{\mathcal{F}}^u_{\sigma}$) coincides with  ${\mathcal{F}}^s_{\mathcal A}$ (${\mathcal{F}}^u_{\mathcal A}$) on $U_\sigma\cap U_{\mathcal A}$, coincides with ${\mathcal{F}}^s_{\mathcal D}$ (${\mathcal{F}}^u_{\mathcal D}$)on  $U_\sigma\cap U_{\mathcal D}$ and coincides with ${\mathcal{F}}^s_{\sigma}$ (${\mathcal{F}}^u_{\sigma}$) out of the set $U_\mathcal A\cup U_{\mathcal D}$ (see figure \ref{conj}). Denote by $\tilde{\pi}^s_{\sigma}: U_{\sigma}\to W^{u}_{\sigma}$ ($\tilde{\pi}^u_{\sigma}: U_{\sigma}\to W^{s}_{\sigma}$) a projection along the leaves of the foliation $\tilde{\mathcal{F}}^s_{\sigma}$ ($\tilde{\mathcal{F}}^u_{\sigma}$). Construct similarly the foliation $\tilde{\mathcal{F}}^{s}_{\sigma'}$ ($\tilde{\mathcal{F}}^{u}_{\sigma'}$) and define the projection  
$\tilde{\pi}^s_{\sigma'}: U_{\sigma'}\to W^{u}_{\sigma'}$ ($\tilde{\pi}^u_{\sigma'}: U_{\sigma'}\to W^{s}_{\sigma'}$) in the neighbourhood $U_{\sigma'}$. Denote by $\tilde{\pi}^s_{\Omega^1},\tilde{\pi}^u_{\Omega^1},\tilde{\pi}^s_{\Omega'^1},\tilde{\pi}^s_{\Omega'^1}$ maps consisting of $\tilde{\pi}^s_{\sigma},\tilde{\pi}^u_{\sigma},\tilde{\pi}^s_{\sigma'},\tilde{\pi}^s_{\sigma'},\sigma\in\Omega^1$, accordingly.

{\bf Step 8.} For each point $a \in\mathcal {A}$ let us define a homeomorphism ${\varphi}_{U_a}:U_a\to U_{a'}$ by the following way. Denote by $U^+_a$ and $U^-_a$ the connected components of $U_a\setminus l_a$ following a rule that any point $z=(z^u_{x},0) \in U_a$ belongs to $U^+_a$ if $z^{u}_{x}>a^{u}_{x}$ and belongs to $U^-_a$ if $z^{u}_{x}<a^{u}_{x}$. Similarly denote the connected components of $U_{a'}\setminus l_{a'}$. Define a homeomorphism ${\varphi}_{U^+_a}:U^+_a\to U^+_{a'}$ by the following way: for a point $z \in U^+_a$ set ${\varphi}_{U^+_a}(z)=z'$, where $z' \in U^+_{a'}$ is the intersection point of the leaves $(\tilde{\pi}^{s}_{\sigma^s_{a'}})^{-1}({\varphi}_{\sigma^s_a}^u(\tilde{\pi}^{s}_{\sigma^s_a}(z)))$ and $(\tilde{\pi}^{u}_{\sigma^u_{a'}})^{-1}({\varphi}_{\sigma^u_a}^s(\tilde{\pi}^{u}_{\sigma^u_a}(z)))$. In the similar way let us define a homeomorphism ${\varphi}_{U^-_a}:U^-_a\to U^-_{a'}$. Set
$${\varphi}_{U_a}(z)=
\begin{cases}{\varphi}_{U^+_a}(z), z\in U^+_a; \cr
{\varphi}_{U^-_a}(z), z\in U^-_a; \cr
{\varphi}_{l_a},z\in l_a.\cr\end{cases}$$

Define a homeomorphism $\varphi_{U_{\mathcal{A}}}:U_{\mathcal{A}} \to U_{\mathcal{A}'}$ as coinciding with ${\varphi}_{U_{a}}$ for each $a \in \mathcal{A}$.

For each point $d\in \mathcal D$ define a homeomorphism ${\varphi}_{U_d}: U_d\to U_{d'}$ by the following way: ${\varphi}_{U_d}(z)$ is the intersection point of the leaves  $(\tilde{\pi}^s_{\sigma^u_{d'}})^{-1}
({\varphi}_{\sigma^u_d}^u(\tilde{\pi}^s_{\sigma^u_d}(z)))$ and
$(\tilde{\pi}^u_{\sigma^u_{d'}})^{-1}({\varphi}_{\sigma^u_d}^s(\tilde{\pi}^u_{\sigma^u_d}(z)))$ belonging to $U_{d'}$. Let us define a homeomorphism  $\varphi_{U_{\mathcal{D}}}:U_{\mathcal{D}} \to U_{\mathcal{D}'}$ as a homeomorphism coinciding with ${\varphi}_{U_{d}}$ for each $d \in \mathcal{D}$.

For $\delta\in\{u,s\}$ denote by $\varphi^\delta_{\Omega^{1\delta}}:W^\delta_{\Omega^{1\delta}}\to W^\delta_{\Omega'^{1\delta}}$ a homeomorphism conjugating the diffeomorphisms $f\vert_{W^\delta_{\Omega^{1\delta}}}$, $f'\vert_{W^\delta_{\Omega'^{1\delta}}}$, coinciding with the homeomorphism $\varphi_{U_{\mathcal A}}$ on $W^\delta_{\Omega^{1\delta}}\cap U_{\mathcal A}$, coinciding with the homeomorphism $\varphi_{U_{\mathcal D}}$ on $W^\delta_{\Omega^{1\delta}}\cap U_{\mathcal D}$ and coinciding with the homeomorphism $\varphi$ out of some neighborhood of the set $W^\delta_{\Omega^{1\delta}}\cap(U_{\mathcal A}\cup U_{\mathcal D})$. Denote by $\varphi^u_{\Omega^1}:W^u_{\Omega^1}\to W^u_{\Omega^1}$ a homeomorphism composed from $\varphi^u_{\Omega^{1u}}$ and $\varphi^u_{\Omega^{1s}}$. Denote by $\varphi^s_{\Omega^1}:W^s_{\Omega^1}\to W^s_{\Omega^1}$ a homeomorphism composed from $\varphi^s_{\Omega^{1u}}$ and $\varphi^s_{\Omega^{1s}}$. 

Set $U_{\Omega^1}=\bigcup\limits_{\sigma\in\Omega^1}U_\sigma$ and $U_{\Omega'^1}=\bigcup\limits_{\sigma'\in\Omega'^1}U_{\sigma'}$. Define a homeomorphism $\varphi_{U_{\Omega^1}}:U_{\Omega^1}\to U_{\Omega'^1}$ as a homeomorphism conjugating the diffeomorphisms $f\vert_{U_{\Omega^1}}$ and $f'\vert_{U_{\Omega'^1}}$, coinciding with the homeomorphism $\varphi_{U_{\mathcal A}}$ on $U_{\Omega^1}\cap U_{\mathcal A}$, coinciding with the homeomorphism $\varphi_{U_{\mathcal D}}$ on $U_{\Omega^1}\cap U_{\mathcal D}$ and such that for a $z\in(U_\sigma\setminus(U_{\mathcal A}\cup U_{\mathcal D}))$, $\varphi_{U_{\Omega^1}}(z)$ is the intersection point of the leaves $(\tilde{\pi}^s_{\Omega^1})^{-1}
({\varphi}_{\Omega^1}^u(\tilde{\pi}^s_{\Omega^1}(z)))$ and
$(\tilde{\pi}^u_{\Omega^1})^{-1}({\varphi}_{\Omega^1}^s(\tilde{\pi}^u_{\Omega^1}(z)))$. 

{\bf Step 9.} For any $t\in(0,1)$ set $U^t_{\mu,\lambda}=\{(x,y)\in\mathbb{R}^2: \vert x\vert {\vert y\vert}^{-\log_{_{\lambda}}\mu} \leq t\}$. For any $\sigma\in\Omega^1$ set  $U^t_\sigma=\psi^{-1}_{\sigma}(U^t_{\mu_\sigma,\lambda_\sigma})$ and $U^t_{\Omega^1}=\bigcup\limits_{\sigma\in\Omega^1}U^t_\sigma$. 

Let us choose a value $t_0\in(0,1)$ such that  
${\varphi}_{U_{\Omega^1}}(U^{t_0}_{\Omega^1})\subset
({\varphi}(U_{\Omega^1})\cup W^s_{\Omega'^{1u}}\cup W^u_{\Omega'^{1s}})$. Set 
${Q}=U_{\Omega^1}\setminus
int~U^{t_0}_{\Omega^1}$, $R=\partial U_{\Omega^1}$, $R_0=\partial U^{t_0}_{\Omega^1}$,  ${Q}'={\varphi}(U_{\Omega^1})\setminus int~{\varphi}_{\Omega^1}(U^{t_0}_{\Omega^1})$,
$R'=\varphi(\partial U_{\Omega^1})$, $R'_0=\varphi_{U_{\Omega^1}}(\partial U^{t_0}_{\Omega^1})$, $\hat Q=p_{_f}(Q)$, $\hat Q'=p_{_{f'}}(Q')$ and  $\hat\varphi_{U_{\Omega^1}}=p_{_{f'}}{\varphi}_{U_{\Omega^1}}
(p_{_f}\vert_{R_0})^{-1}:\hat R_0\to \hat R'_0$. 

By the construction the sets $\hat Q$, $\hat Q'$ have the same number of the connected components each of them is homeomorphic to the standard two-dimensional annulus 
(see figure \ref{k}, where the set $\hat Q$ is coloured). Then there is a homeomorphism  $\hat\varphi_{\hat Q}:\hat Q\to \hat Q'$ such that  $\hat\varphi_{\hat Q}\vert_{\hat R}=\hat\varphi$ and $\hat\varphi_{\hat Q}\vert_{\hat R_0}=\hat\varphi_{U_{\Omega^1}}$.

\begin{figure} \begin{center} \epsfig
{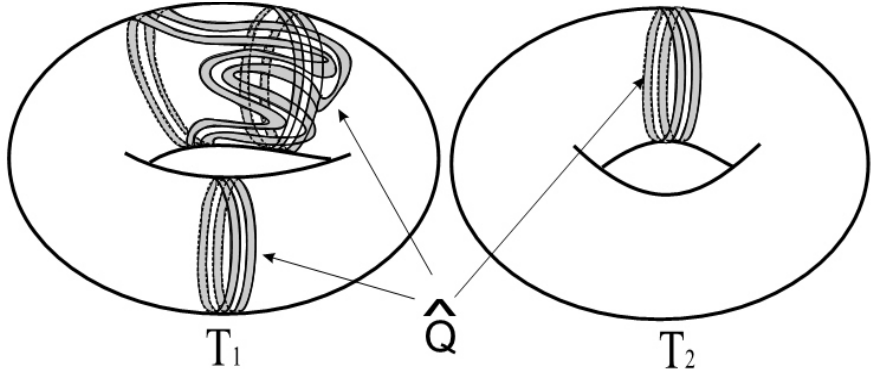} \caption{Illustration to Step 9} \label{k} \end{center} 
\end{figure}

Denote by ${\varphi}_Q:Q\to {Q}'$ a lift of the homeomorphism $\hat\varphi_{\hat Q}$ coinciding with ${\varphi}$ on $\partial U_{\Omega^1}$. Define a homeomorphism  $\varphi_3:V_{f} \to V_{f'}$ by the  formula:
$$\varphi_3(x)=\begin{cases}\varphi_{U_{\Omega^1}}(x), x\in
U^{t_0}_{\Omega^1}; \cr {\varphi}_{Q}(x),x\in
{Q}; \cr {\varphi}(x), x\in
M^2\setminus {U}_{\Omega^1}.\cr\end{cases}$$
Let us define a homeomorphism $h:M^2\setminus(\Omega^0\cup\Omega^2)\to M'^2\setminus(\Omega'^0\cup\Omega'^2)$  by the formula:
$$h(x)=\begin{cases}\varphi_{2}(x), x\in
(W^s_{\mathcal L^s}\setminus\mathcal L^s \cup W^u_{\mathcal L^u}\setminus\mathcal L^u); 
\cr {\varphi}_{3}(x),x\in
(V_f\setminus (W^s_{\mathcal L^s}\cup W^u_{\mathcal L^u}); 
\cr {\varphi}_{\mathcal L^s}(x), x\in
\mathcal L^s;
\cr {\varphi}_{\mathcal L^u}(x), x\in
\mathcal L^u;
\cr {\varphi}^u_{\Omega^{1s}}(x), x\in
W^u_{\Omega^{1s}};
\cr {\varphi}^s_{\Omega^{1u}}(x), x\in
W^s_{\Omega^{1u}}.\cr\end{cases}$$

Then, to obtain a desired homeomorphism, it suffices to extend the homeomorphism $h$ continuously to the set $\Omega^0\cup\Omega^2$.
\end{demo}


\begin{thebibliography}{99}

\bibitem{ArGr90} S. Aranson, V. Grines.  {\em The topological classification of cascades on closed two-dimensional manifolds}. Russian Mathematical Surveys. 45:1, (1990), 1-35.

\bibitem{BoLa} C. Bonatti, R. Langevin. {\em Diff\'eomorphismes de Smale des
surfaces}. Ast\'erisque. Soci\'et\'e mathematique de France. Paris. 1998, №
250.

\bibitem{BoGrLa01} C. Bonatti, V. Grines and R. Langevin. {\em Dynamical systems in dimension 2 and 3: Conjugacy invariants and classification}. Computational and Applied Mathematics. {\bf  20}, 2001, no. 1-2, 11--50. 

\bibitem{Bow71} R. Bowen. {\em Periodic points and measures for axiom A diffeomorphisms}. Transactions of the American. Math. Soc. 1971. {\bf  154},   337--397.

\bibitem{Gr74} V. Grines.  {\em The topological equivalence of one-dimensional basic sets of diffeomorphisms on two-dimensional manifolds}. Uspekhi Mat. Nauk, {\bf 29}, no 6(180) (1974), 163--164.

\bibitem{Gr75} V. Grines.  {\em The topological conjugacy of diffeomorphisms of a two-dimensional manifold on one-dimensional orientable basic sets}. I. Tr. Mosk. Mat. Obs., 32, MSU, M., 1975, 35--60.

\bibitem{Gr97} V. Grines.  {\em On the topological classification of structurally stable diffeomorphisms of surfaces with one-dimensional attractors and repellers}. Sb. Math. 188, no. 4, (1997), 537-569.

\bibitem{Gr1999} V.Z. Grines, {\em Topological classification of one-dimensional attractors and repellers of A-diffeomorphisms of surfaces by means of automorphisms of fundamental groups of supports}. Dynamical systems. 7. J. Math. Sci. (New York) 95 (1999), no. 5, 2523--€"2545.

\bibitem{GrZhu2006} Grines, V. Z.; Zhuzhoma, E. V. Expanding attractors. Regul. Chaotic Dyn. 11 (2006), no. 2, 225-246.

\bibitem{Ko} C. Kosniowski. {\em A first course in algebraic topology}. CUP Archive, 1980.

\bibitem{MeSt1987} W. de Melo and S. van Strien, {\em Diffeomorphisms on surfaces with a finite number of moduli}. Ergodic Theory Dynam. Systems {\bf 7} (1987), no. 3, 415--462.

\bibitem{MiPo2010} T.M. Mitryakova, O.V. Pochinka, {\em On necessary and sufficient conditions for the topological conjugacy of surface diffeomorphisms with a finite number of orbits of heteroclinic tangency}. (Russian) Tr. Mat. Inst. Steklova 270 (2010), Differentsialnye Uravneniya i Dinamicheskie Sistemy, 198--219; translation in Proc. Steklov Inst. Math. 270 (2010), no. 1, 194--215.

\bibitem{NP73} S. Newhouse and J. Palis, {\em Bifurcations of Morse-Smale dynamical systems}.  Dynamical systems (edited by M. Peixoto). Academic Press, New York, 1973.

\bibitem{Pa78} J. Palis.  {\em A differentiable invariant of topological conjugacy and moduli of
stability}. Asterisque. 1978. {\bf 51}, 335-346.

\bibitem{Plyk74}  R.V. Plykin. {\em  Sources and sinks of A-diffeomorphisms of surfaces}. (Russian) Mat. Sb. (N.S.) {\bf 94}(136) (1974), 243--264, 336.

\bibitem{ChShShTu} L. P. Shilnikov, A. L. Shilnikov, D. V. Turaev, and L. O. Chua, {\em Methods of Qualitative Theory in Nonlinear Dynamics. I}.  (World Sci., Singapore, 1998; Inst. Kompyut. Issled., Moscow, 2004).

\bibitem{Smale67}  S. Smale. {\em Differentiable dynamical systems}. Bull. Amer. Math. Soc. 1967. {\bf 73}, no 1, 741--817.

\label{last}
\end{thebibliography}
\end{document}